\documentclass{amsart}

\usepackage[latin1]{inputenc}
\usepackage{color}
\usepackage{url}
\usepackage{pstricks}
\usepackage{amssymb}
\def\Dfn#1{{\sf #1}}

\def\lcm#1{\operatorname{lcm}\{ #1 \}}

\usepackage{amssymb}
\usepackage{subfigure}
\usepackage{setspace}
\usepackage{graphicx}

\def\Cat{{\rm Cat}}
\def\supp{{\rm supp}}
\def\An{A_{n-1}}
\def\NN{{\rm NN}}
\def\NC{{\rm NC}}
\def\K{{\sf Krew}}
\def\P{{\sf Pan}}
\def\c{{\sf conj}}
\def\T{\mathcal{T}}

\def\O{\mathcal{O}}

\def\id{\mathbf{1}}

\newcommand{\Pan}{\P}
\newcommand{\rot}{\K}
\newcommand{\Ind}{\operatorname{Ind}}

\newtheorem{theorem}{Theorem}[section]
\newtheorem*{maintheorem}{Main Theorem}
\newtheorem*{panyushevconjecture}{Panyushev Conjectures}
\newtheorem{proposition}[theorem]{Proposition}
\newtheorem{corollary}[theorem]{Corollary}
\newtheorem{lemma}[theorem]{Lemma}

\theoremstyle{definition}
\newtheorem{definition}[theorem]{Definition}
\newtheorem{example}[theorem]{Example}

\newtheorem{remark}{Remark}

\begin{document}
  \title[A uniform bijection between $\NN$ and $\NC$.]{A uniform bijection between nonnesting and noncrossing partitions}

  \author{Drew Armstrong}
  \address{Department of Mathematics, University of Miami, Coral Gables, FL, 33146}
  \email{armstrong@math.umiami.edu}

  \author{Christian Stump}
  \address{LaCIM, Universit\'e du Qu\'ebec \`a Montr\'eal, Montr\'eal (Québec), Canada}
  \email{christian.stump@univie.ac.at}

  \author{Hugh Thomas}
  \address{Department of Mathematics and Statistics, University of New Brunswick, Fredericton NB, E3B 5A3}
  \email{hthomas@unb.ca}

  \date{\today}

  \begin{abstract}
    In 2007, D.I.~Panyushev defined a remarkable map on the set of nonnesting partitions (antichains in the root poset of a finite Weyl group). In this paper we identify Panyushev's map with the Kreweras complement on the set of noncrossing partitions, and hence construct the first uniform bijection between nonnesting and noncrossing partitions. Unfortunately, the proof that our construction is well-defined is case-by-case, using a computer in the exceptional types. Fortunately, the proof involves new and interesting combinatorics in the classical types. As consequences, we prove several conjectural properties of the Panyushev map, and we prove two cyclic sieving phenomena conjectured by D.~Bessis and V.~Reiner.
  \end{abstract}

  \maketitle

  \section{Introduction}

  To begin we will describe the genesis of the paper.

  \subsection{Panyushev complementation}
  Let $\Delta\subseteq\Phi^+\subseteq\Phi$ be a triple of {\sf simple roots}, {\sf positive roots}, and a {\sf crystallographic root system} corresponding to a {\sf finite Weyl group} $W$ of \Dfn{rank} $r$. We think of $\Phi^+$ as a poset in the usual way, by setting $\alpha\leq\beta$ whenever $\beta-\alpha$ is in the nonnegative span of the simple roots $\Delta$. This is called the {\sf root poset}. The set of \Dfn{nonnesting partitions} $\NN(W)$ is defined to be the set of antichains (sets of pairwise-incomparable elements) in $\Phi^+$. This name is based on a pictorial presentation of antichains in the classical types. It is well known that the number of nonnesting partitions is equal to the \Dfn{Catalan number}
  \begin{equation*}
    \Cat(W):=\prod_{i=1}^r \frac{d_i+h}{d_i},
  \end{equation*}
  where $d_1\leq d_2\leq \cdots \leq d_r = h$ are the degrees of a fundamental system of polynomial invariants for $W$ (called the {\sf degrees} of $W$), and where $h$ is the \Dfn{Coxeter number}. This formula was first conjectured by Postnikov~\cite[Remark~2]{Rei1997} and at least two uniform proofs are known \cite{Ath2004,Hai1994}. These enumerations were established in somewhat different contexts; the link to the combinatorics of antichains is supplied in both cases by \cite{CP2002}.

  In 2007, Panyushev defined a remarkable map on nonnesting partitions~\cite{Pan2008}. To describe it, we first note than an antichain $I\subseteq\Phi^+$ corresponds bijectively to the {\sf order ideal} $\langle I\rangle\subseteq\Phi^+$ that it generates. The \Dfn{Panyushev complement} is defined as follows.

\begin{definition}
Given an antichain of positive roots $I\subseteq\Phi^+$, define ${\sf Pan}(I)$ to be the antichain of minimal roots in $\Phi^+\setminus\langle I\rangle$.

\end{definition} 

\begin{figure}
\begin{center}
\includegraphics[scale=.5]{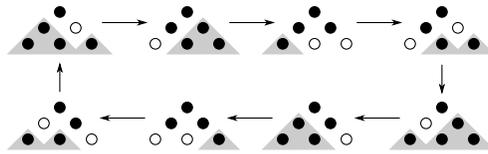}
\end{center}
\caption{An orbit of the Panyushev complement}
\label{fig:panorbit}
\end{figure}

For example, Figure \ref{fig:panorbit} displays a single orbit of the Panyushev complement acting on the root poset of type $A_3$.  The antichain in each picture corresponds to the maximal black dots in the order ideal given by the shaded area. In \cite[Conjecture~2.1]{Pan2008} Panyushev made several conjectures about the Panyushev complementation, which have remained open even in type $A$. Even though $\P$ can be defined on any poset, those conjectures provide strong evidence that the Panyushev complementation behaves in a very special way for root posets, and that it has a particular meaning in this case which has not yet been explained.

\begin{panyushevconjecture}
	Let $W$ be a finite Weyl group of rank $r$, with $h$ its Coxeter number, and $\P$ the Panyushev complement on antichains in the associated root poset $\Phi^+$. Moreover, let $\omega_0$ be the unique longest element in $W$.
	\begin{itemize}
		\item[(i)] $\P^{2h}$ is the identity map on $\NN(W)$,
		\item[(ii)] $\P^h$ acts on $\NN(W)$ by the involution induced by $-\omega_0$,
		\item[(iii)] For any orbit $\O$ of the Panyushev complement acting on $\NN(W)$, we have
		$$\frac{1}{|\O|}\sum_{I \in \O}|I| = r/2.$$
	\end{itemize}
\end{panyushevconjecture}

For example, in type $A_3$ we have $2h=8$, and the Panyushev complement has three orbits, of sizes $2$, $4$, and $8$ (the one pictured). In type $A$, $\omega_0$ acts by $\alpha_i \mapsto -\alpha_{n-i}$ where $\alpha_i$ denotes the $i$-th simple root in the linear ordering of the Dynkin diagram. It can be easily seen in the pictured orbit, that $\P^h$ acts by ``flipping'' the root poset (this corresponds to reversing the linear ordering of the Dynkin diagram), and 
that $\P^{2h}$ is the identity map. Moreover, the average number of elements in this orbit is $\frac{1}{8}(2+1+1+2+2+1+1+2) = 3/2$.

In this paper we will prove the following.
\begin{theorem}
  The Panyushev Conjectures are true.
\end{theorem}
However, the proof of this theorem is not the main goal of the paper. Instead, we will use the Panyushev complement as inspiration to solve an earlier open problem: to find a uniform bijection between the antichains in $\Phi^+$ and a different sort of Catalan object, the {\sf noncrossing partitions}. We will then use the combinatorics we have developed to prove the Panyushev Conjectures.

\subsection{Kreweras complementation} There is also a notion of {\sf noncrossing partitions} for root systems, which we now describe.

Let $T$ be the set of all \Dfn{reflections} in a \Dfn{finite Coxeter group} $W$. Those are given by the reflections defined by the positive roots in a (not necessarily crystallographic) finite root system $\Phi$. Let $c \in W$ be a \Dfn{Coxeter element} (i.e., the product of the \Dfn{simple reflections} $S$ defined by the simple roots in some order). Then the set of {\sf noncrossing partitions} is
\begin{equation*}
  \NC(W,c):=\{ w\in W: \ell_T(w)+\ell_T(cw^{-1})=r\}\subseteq W,
\end{equation*}
where $r$ is the rank of $W$. For a full exposition of this object and its history, see~\cite{Arm2006}. It turns out that $\NC(W,c)$ is also counted by the Catalan number $\Cat(W)$, but in this case {\bf no uniform proof is known} (the only proof is case-by-case, using a computer for the exceptional types). In this paper we will (partially) remedy the situation by constructing a uniform bijection between antichains in $\Phi^+$ and the noncrossing partitions $\NC(W,c)$. It is only a partial remedy because our proof that the construction is well-defined remains case-by-case.

Our bijection relies on the Panyushev complement and a certain map on noncrossing partitions, which we now describe. The type $A$ noncrossing partitions were first studied in detail by Kreweras~\cite{Kre1972}, as pictures of ``noncrossing partitions'' of vertices around a circle. He noticed that the planarity of these pictures yields a natural automorphism, which we call the {\sf Kreweras complement}.

\begin{definition}
Given a noncrossing partition $w\in\NC(W,c)\subseteq W$, let ${\sf Krew}(w):=cw^{-1}$. Since the reflection length $\ell_T$ is invariant under conjugation it follows that ${\sf Krew}(w)$ is also in $\NC(W,c)$.
\end{definition}

\begin{figure}
\begin{center}
\includegraphics[scale=.5]{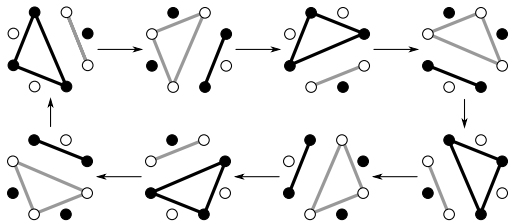}
\end{center}
\caption{An orbit of the Kreweras map}
\label{fig:kreworbit}
\end{figure}

In type $A_{n-1}$, the set $\NC(W,c)$ consists of partitions of the vertices $\{1,2,\ldots,n\}$ placed around a circle, such that the convex hulls of its equivalence classes are nonintersecting (``noncrossing''). To describe the classical Kreweras map, we place vertices $\{1',1,2',2,\ldots,n',n\}$ around a circle; if $\pi$ is a noncrossing partition of $\{1,2,\ldots,n\}$ then ${\sf Krew}(\pi)$ is defined to be the {\bf coarsest} partition of $\{1',2',\ldots,n'\}$ such that $\pi\cup{\sf Krew}(\pi)$ is noncrossing. For example, Figure \ref{fig:kreworbit} shows a single orbit of ${\sf Krew}$ acting on the noncrossing partitions of a square (given by the black vertices). Note here that ${\sf Krew}^2$ rotates the square by $90^\circ$.

For a general root system we have ${\sf Krew}^2(w)=cwc^{-1}$; that is, ${\sf Krew}^2$ is conjugation by the Coxeter element. Since any Coxeter element $c$ has order $h$ (indeed this is an equivalent definition of the Coxeter number $h$) we conclude that ${\sf Krew}^{2h}$ is the identity map. Thus we prove part (i) of the Panyushev conjectures by finding a bijection from antichains to noncrossing partitions that sends ${\sf Pan}$ to ${\sf Krew}$.

\subsection{Panyushev complementation $=$ Kreweras complementation} Since no uniform bijection currently exists, we will create one, essentially by {\em declaring} that ${\sf Pan}={\sf Krew}$. The key observation is the following.

Since a Dynkin diagram of finite type is a tree, we may partition the simple reflections $S$ into sets $S=L\sqcup R$ such that the elements of $L$ commute pairwise, as do the elements of $R$. Let $c_L$ denote the product of the reflections $L$ (in any order) and similarly let $c_R$ denote the product of the reflections $R$. Thus, $c_L$ and $c_R$ are involutions in $W$ and $c=c_Lc_R$ is a special Coxeter element, called a {\sf bipartite Coxeter element}.

The data for ${\sf Pan}$ consists of a choice of simple system $\Delta$, which from now on we will partition as $\Delta=\Delta_L\sqcup\Delta_R$; and the data for ${\sf Krew}$ consists of a Coxeter element, which from now on we will assume to be $c=c_Lc_R$. With this in mind, Panyushev observed that his map has two distinguished orbits: one of size $h$ which consists of the sets of roots at each rank of the root poset; and one of size 2, namely $\{\Delta_L,\Delta_R\}$. Similarly, the Kreweras map on $\NC(W,c_Lc_R)$ has two distinguished orbits: one of size $h$ consisting of
\begin{equation*}
c_L, c_Lc_Rc_L, \ldots, c_Rc_Lc_R, c_R;
\end{equation*}
and one of size 2, namely $\{{\bf 1},c\}$. The attempt to match these orbits was the genesis of our Main Theorem.

To understand its statement, we must first discuss parabolic recursion. Let $W_J\subseteq W$ denote the {\sf parabolic subgroup} generated by some subset $J\subseteq S$ of simple reflections, and let $\Delta_J\subseteq\Phi_J^+\subseteq\Phi^+$ be the corresponding simple and positive roots. Antichains and noncrossing partitions may be restricted to $W_J$ as follows. Let $I\subseteq\Phi^+$ be an antichain and define its support ${\sf supp}(I)=\langle I\rangle\cap\Delta$ to be the simple roots below it. If ${\sf supp}(I)\subseteq J$ then $I$ is also an antichain in the parabolic sub-root system $\Phi^+_J$. Similarly, the set $J$ induces a unique partition of the diagram $J=L_J\sqcup R_J$ with $L_J\subseteq L$ and $R_J\subseteq R$, and we may discuss the {\sf parabolic noncrossing partitions}
\begin{equation*}
  \NC(W_J,c_{L_J}c_{R_J})\subseteq\NC(W,c_Lc_R).
\end{equation*}
With these notions in mind, we state our main theorem.

\begin{maintheorem}\label{th:main}
Let $S=L\sqcup R$ be a bipartition of the simple reflections with corresponding bipartition $\Delta=\Delta_L\sqcup\Delta_R$ of the simple roots and bipartite Coxeter element $c=c_Lc_R$. Then {\bf there exists a (unique) bijection} $\alpha$ from nonnesting partitions $\NN(W)$ to noncrossing partitions $\NC(W,c)$ satisfying the following three properties:
\begin{itemize}
\item $\alpha(\Delta_L)={\bf 1}$,\hfill {\em (initial condition)}
\item $\alpha\circ{\sf Pan}={\sf Krew}\circ\alpha$, \hfill {\em (${\sf Pan}={\sf Krew}$)}
\item $\alpha(I)=\left(\prod_{s\in L\setminus{\sf supp}(I)} s\right) \alpha |_{{\sf supp}(I)}(I).$ \hfill {\em (parabolic recursion)}
\end{itemize}
\end{maintheorem}

That is, to compute $\alpha$ of an antichain $I$, let $J={\sf supp}(I)$. If $J\subsetneq S$ then we think of $I$ as an antichain in the proper subsystem $\Phi^+_J$. We compute $\alpha|_J(I)$, which is an element of
\begin{equation*}
\NC(W_J,c_{L_J}c_{R_J})\subseteq W_J,
\end{equation*}
and then multiply on the left by the simple $L$-reflections {\bf not} in $J$. As $J$ gets smaller, we will reach the initial condition $\alpha(\Delta_{L_J})={\bf 1}$. If $J=S$ then we apply the map ${\sf Pan}$ $k$ times until we have ${\sf supp}({\sf Pan}^k I)\subsetneq S$. Finally, we apply $\alpha$ and then ${\sf Krew}^{-k}$. That this process works is the content of the theorem.

\medskip

\begin{remark}
  The statement of the Main Theorem is {\bf uniform}. (That is, it is expressed purely in terms of root systems.) Unfortunately, we will prove the theorem in a case-by-case way. Fortunately, the proof involves new and interesting combinatorics in the classical types. (Which is new and interesting even in type $A$.)
\end{remark}
\medskip

We note that the interaction between the ``nonnesting'' and ``noncrossing'' properties is a subtle phenomenon, even in type $\An$ alone (see \cite{CDDSY2007}). There has also been earlier progress on the problem for general finite root systems: A.~Fink and B.I.~Giraldo \cite{FG2009} and M.~Rubey and the second author~\cite{RS2010} have both constructed bijections which work for the classical types. These bijections have an advantage over ours in that they both preserve the ``parabolic type'' of noncrossing and nonnesting partitions. However, our bijection has the advantage of being uniform for root systems, as well as proving the Panyushev conjectures and a cyclic sieving phenomenon as described in the following section.

  \subsection{Cyclic Sieving} The cyclic sieving phenomenon was introduced by V.~Reiner, D.~Stanton, and D.~White in \cite{RSW2004} as follows: let $X$ be a finite set, let $X(q) \in \mathbb{Z}[q]$ and let $\mathcal{C}_d = \langle c \rangle$ be a cyclic group of order $d$ acting on X. The triple $(X,X(q),\mathcal{C}_d)$ exhibits the \Dfn{cyclic sieving phenomenon (CSP)} if
  $$ [X(q)]_{q = \zeta^k} = \big| X^{c^k} \big|,$$
  where $\zeta$ denotes a primitive $d$-th root of unity and $X^{c^k} := \{x \in X : c^k(x) = x\}$ is the fixed-point set of $c^k$ in $X$. Let 
  \begin{align}
    X(q) &\equiv a_0 + a_1 q + \ldots + a_{d-1} q^{d-1} \mod (q^d-1). \label{eq:CSPorbitlengths}
  \end{align}
  An equivalent way to define the CSP is to say that $a_i$ equals the number of $\mathcal{C}_d$-orbits in $X$ whose stabilizer order divides $i$ \cite[Proposition 2.1]{RSW2004}.

  Bessis and Reiner recently showed that the action of the Coxeter element on noncrossing partitions together with a remarkable $q$-extension of the Catalan numbers $\Cat(W)$ exhibits the CSP: define the \Dfn{$q$-Catalan number}
  \begin{equation*}
    \Cat(W;q):=\prod_{i=1}^r \frac{[d_i+h]_q}{[d_i]_q},
  \end{equation*}
  where $[k]_q=1+q+q^2+\cdots +q^{k-1}$ is the usual {\sf $q$-integer}. It is not obvious, but it turns out (see Berest, Etingof, and Ginzburg \cite{BEG2003}) that this number is a polynomial in $q$ with nonnegative coefficients. In type $A_{n-1}$, the formula reduces to the classical $q$-Catalan number of F\"urlinger and Hofbauer \cite{FH1985}. That is, we have 
  \begin{equation*}
    \Cat(A_{n-1};q)=\frac{1}{[n+1]_q}\begin{bmatrix} 2n\\ n \end{bmatrix}_q,
  \end{equation*}
  where $\left[\begin{smallmatrix} a \\ b \end{smallmatrix}\right]_q=\frac{[a]_q!}{[b]_q![a-b]_q!}$ is the {\sf Gaussian binomial coefficient} and $[k]_q!:=[1]_q[2]_q\cdots [k]_q$ is the {\sf $q$-factorial}.

  For a Coxeter element $c\in W$, it follows directly from the definition that the map $\c(w) = cwc^{-1}$ is a permutation of the set $\NC(W,c)$ of noncrossing partitions. In classical types, this corresponds to a ``rotation'' of the pictorial presentation.

  \begin{theorem}[Bessis and Reiner \cite{BR2007}] \label{th:BR}
    The triple $\big(\NC(W),\Cat(W;q), \langle\c\rangle\big)$ exhibits the CSP for any finite Coxeter group $W$.
  \end{theorem}

  Actually, they proved this result in the greater generality of {\sf finite complex reflection groups}; we will restrict the current discussion to {\sf (crystallographic) finite real reflection groups} --- that is, finite Coxeter groups and finite Weyl groups, respectively. At the end of their paper, Bessis and Reiner \cite{BR2007} conjectured several other examples of cyclic sieving, two of which we will prove in this paper.

  \begin{theorem}
    Let $W$ be a finite Coxeter group respectively finite Weyl group.
    \begin{itemize}
      \item[(i)] The triple $\big(\NC(W),\Cat(W;q), \langle \K \rangle\big)$ exhibits the CSP.
      \item[(ii)] The triple $\big(\NN(W),\Cat(W;q), \langle \P \rangle\big)$ exhibits the CSP.
    \end{itemize}
    \label{th:CSPNC}
  \end{theorem}

  Note that (i) is a generalization of Theorem \ref{th:BR} since ${\sf Krew}^2$ is the same as conjugation by the Coxeter element. The type $A$ version of (i) has been proved by D.~White (see \cite{BR2007}) and independently by C.~Heitsch \cite{Hei}; C.~Krattenthaler has announced a proof of a more general version for complex reflection groups which appeared in the exceptional types in \cite{KM2008}; and will appear for the group $G(r,p,n)$ in \cite{kra}. In this paper we find it convenient to present an independent proof, on the way to proving our Main Theorem. Combining (i) and the Main Theorem then yields (ii) as a corollary.

  \subsection{Outline} The paper is organized as follows.

  In {\bf Section 2}, we introduce a notion of {\sf noncrossing handshake configurations} for the classical types, and define a bijection $\phi_W$ from noncrossing handshake configurations $\T_W$ to noncrossing partitions $\NC(W,c)$.  We establish the cyclic sieving phenomenon for noncrossing partitions using these bijections in classical type, and via a computer check for the exceptional types.

  In {\bf Section 3}, we define a bijection $\psi_W$ from the nonnesting partitions of $W$ to $\T_W$ in the classical types. Using this, we establish the cyclic sieving phenomenon for nonnesting partitions in the classical types, and again via a computer check for the exceptional types.

  In {\bf Section 4}, we show that the bijection from the nonnesting partitions of $W$ to $\T_W$ in the classical types satisfies a suitable notion of parabolic induction.

  In {\bf Section 5}, we put together the bijections from sections two and three to prove the Main Theorem.  The calculations for the exceptional types were done using \texttt{Maple} code, which is available from the first author.

	In the final section, {\bf Section 6}, we use the combinatorics describing the Panyushev and the Kreweras complementation to prove the {\bf Panyushev conjectures}.

\section{The Kreweras CSP for noncrossing partitions}\label{sectiontwo}

  In this section, we prove Theorem~\ref{th:CSPNC}(i) for every type individually. For type $\An$, C.~Heitsch proved the theorem by connecting noncrossing partitions of type $\An$ to noncrossing set partitions of $[n] := \{1,\ldots,n\}$ and moreover to noncrossing handshake configurations of $[2n]$ and to rooted plane trees. For the classical types, we will explore a connection which is related to the construction of C.~Heitsch as described in Remark~\ref{re:heitschbijection}.
  
\subsection{Type $A$}
  
\begin{figure}
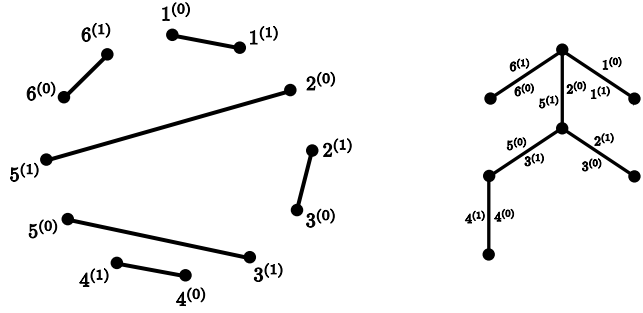

    \centering
    \begin{tabular}{ccc}
      \input{handshake6.tex}
      & \hspace{20pt} &
      \input{rootedplanetree6.tex} 
    \end{tabular}
    \caption{The noncrossing handshake configuration $T \in \T_6$ for $w = (2,3,5)$ and its associated rooted plane tree.}
    \label{fig:handshakeA}
  \end{figure}

  Fix the linear Coxeter element $c$ to be the long cycle $(1,2,\ldots,n)$. Here, \emph{linear} refers to the fact that it comes from a linear ordering of the Dynkin diagram. It is well-known that the set of noncrossing partitions $\NC_n := \NC(\An)$ can be identified with the set of noncrossing handshake configurations. The ground set consists of $2$ copies of $[n]$ colored by $0$ and $1$ drawn on a circle in the order $1^{(0)},1^{(1)},\ldots,n^{(0)},n^{(1)}$. A \Dfn{noncrossing handshake configuration} is defined to be a noncrossing matching of those $2$ copies of $[n]$, see Figure~\ref{fig:handshakeA}. As shown in the figure, they are in natural bijection with rooted plane trees.
  The bijection $\phi_{\An} : \T_n \longrightarrow \NC_n$ is then, for $w = \phi_{\An}(T)$, given by
  $$ \big(i^{(1)},j^{(0)}\big) \in T \Leftrightarrow w(i) = j.$$
  For a direct description of noncrossing partitions in terms of rooted plane trees see e.g. \cite[Figure 6]{Ber2007}.
  \begin{remark}\label{re:CoxeterElements}
    Observe that the described construction does not require the choice of the linear Coxeter element. As the Coxeter elements in type $\An$ are exactly the long cycles, one obtains analogous constructions by labelling the vertices of $\T_n$ by any given long cycle. This corresponds to the natural isomorphism between $\NC(W,c)$ and $\NC(W,c')$ given by conjugation sending $c$ to the Coxeter element $c'$. We will make use of this flexibility later on in this paper.
  \end{remark}
  The following proposition follows immediately from the definition.
  \begin{proposition} \label{prop:KrewerasA}
    The Kreweras complementation on $\NC_n$ can be described in terms of $\T_n$ by clockwise rotation of all edges by one, or, equivalently, by counterclockwise rotation of all vertex labels by one. I.e., for $T \in \T_n$, we have
    $$
      \big(i^{(1)},j^{(0)}\big) \in T \Leftrightarrow \big(j^{(1)},(i+1)^{(0)}\big) \in \K(T).
    $$
  \end{proposition}

  \begin{remark}\label{re:heitschbijection}
    One can easily deduce the proposition as well from O.~Bernardi's description \cite[Figure~6]{Ber2007} and the definition of the Kreweras complementation of a set partition to be its coarsest complementary set partition. C.~Heitsch obtains analogous results in \cite{Hei} by directly considering a bijection $\phi'$ between $\T_n$ and $\NC_n$ which is related to the bijection $\phi$ described above by $\phi'(w) = \phi(\K(w))$.
  \end{remark}

  For more readability, we set $\Cat_n(q) := \Cat(\An;q)$, and $\Cat_n := \Cat_n(1)$.
  \begin{theorem} \label{th:CSPNCA}
    The triple $\big(\NC_n,\Cat_n(q),\langle \K \rangle \big)$ exhibits the CSP.
  \end{theorem}
  \begin{proof}
    The theorem follows immediately from \cite[Theorem 8]{Hei2}: let $d$ be an integer such that $d \big| 2n$ and let $\zeta$ be a primitive $d$-th root of unity. Then it follows e.g. from \cite[Lemma 3.2]{EF2008} that $\Cat_n(q)$ reduces for $q = \zeta$ to
    \begin{eqnarray}
      \big[\Cat_n(q)\big]_{q=\zeta} =
        \left\{ \begin{array}{cl}
          {\displaystyle\Cat_n} &\mbox{if $d=1$} \\[6pt]
          n {\displaystyle \Cat_{\frac{n-1}{2}}} & \mbox{if $d=2$ and $n$ odd} \\[6pt]
          {\displaystyle \binom{2n/d}{n/d}} &\mbox{if $d \geq 2$, $d\big|n$} \\[6pt]
          {\displaystyle 0} & \mbox{otherwise}
        \end{array} \right. \label{eq:CatA}
    \end{eqnarray}
    In \cite[Theorem 8]{Hei2}, C.~Heitsch proved that noncrossing handshake configurations of $2n$ which are invariant under a $d$-fold rotation, i.e., for which $\K^{2n/d}(T) = T$, are counted by those numbers.
  \end{proof}

\subsection{Types $B$ and $C$}

  As the reflection groups of types $B$ and $C$ coincide, the notions of noncrossing partitions do as well. Therefore we restrict our attention to type $C$. In this case, we fix the linear Coxeter element $c$ to be the long cycle $(1,\ldots,n,-1,\ldots,-n)$ and keep in mind that we could replace $c$ by any long cycle of analogous form. $\NC(C_n)$ can be seen as the subset of $\NC(A_{2n-1})$ containing all elements for which $i \mapsto j$ if and only if $-i \mapsto -j$, where $n+i$ and $-i$ are identified. 
$\T_{C_n}$ is 
defined to be the set of all noncrossing handshake configuration $T$ of $[\pm n]$ for which $(i^{(1)},j^{(0)}) \in T$ if and only if $(-i^{(1)},-j^{(0)}) \in T$. The Kreweras complementation on $\NC(C_n)$ is again the clockwise rotation of all edges by $1$. Observe that the symmetry property is expressed in terms of the Kreweras complementation by $\K^{2n}(T) = T$ for $T \in \T_{C_n}$. In particular, we see that the Kreweras map of order $4n$ on $\T_{C_n}$ is never free. By construction, the bijection $\phi_{A_{2n-1}} : \T_{2n} \tilde{\longrightarrow} \NC_{2n}$ restricts to a bijection
  $$\phi_{C_n} : \T(C_n) \tilde{\longrightarrow} \NC(C_n),$$
  which is compatible with the Kreweras complementation, i.e.,
  $$\phi_{C_n}(\K(T)) = \K(\phi_{C_n}(T)).$$
  
  For the proof of Theorem~\ref{th:CSPNC}(i) in type $C$, we need the following observation.
  \begin{lemma}
    Let $d_1,d_2 \big| 2n$ and let $d_3 = \lcm{d_1,d_2}$. $T \in \T_n$ is invariant both under $d_1$- and $d_2$-fold rotation if and only if $T$ is invariant under $d_3$-fold symmetry.
  \end{lemma}
  
  \begin{proposition}
    The triple $\big(\NC(C_n),\Cat(C_n;q),\langle \K \rangle \big)$ exhibits the CSP.
  \end{proposition}
  The proof in type $C$ is a simple corollary of the proof in type $A$.
  \begin{proof}
    The $q$-Catalan number $\Cat(W;q)$ reduces for $W = C_n$ to
    \begin{align*}
      \Cat(C_n,q) &= \begin{bmatrix} 2n\\ n \end{bmatrix}_q.
    \end{align*}
    Let $d$ be an integer such that $d \big| 4n$ and let $\zeta$ be a primitive $d$-th root of unity. Then it follows again from \cite[Lemma 3.2]{EF2008} that $\Cat(C_n,q)$ reduces for $q = \zeta$ to
    \begin{eqnarray*}
      \big[\Cat(C_n,q)\big]_{q=\zeta} = \left\{ \begin{array}{cl} {\displaystyle \binom{4n/d}{2n/d}} &\mbox{if $d$ even and $d\big|2n$} \\[10pt] {\displaystyle \binom{2n/d}{n/d}} &\mbox{if $d$ odd} \\[6pt] {\displaystyle 0} & \mbox{otherwise} \end{array} \right.
    \end{eqnarray*}
    Let $d\big|4n$. Then by the previous lemma, the number of elements in $\T_{C_n}$ which are invariant under $d$-fold symmetry, i.e., for which $\K^{4n/d}(T) = T$, are exactly those elements in $\T_{2n}$ which are invariant under $\lcm{d,2}$-fold symmetry. The proposition follows.
  \end{proof}

\subsection{Type $D$}
  In this case, we fix the linear Coxeter element $c$ to be $(1,\ldots,n-1,-1,\ldots,-n+1)(n,-n)$. As in types $A$ and $C$, the noncrossing handshake configuration in type $D$ comes from noncrossing set partitions of type $D$ as defined in \cite{AR2004} by replacing every point $i$ by the two points $i^{(0)}$ and $i^{(1)}$, together with the appropriate restrictions, as described below.

\begin{figure}
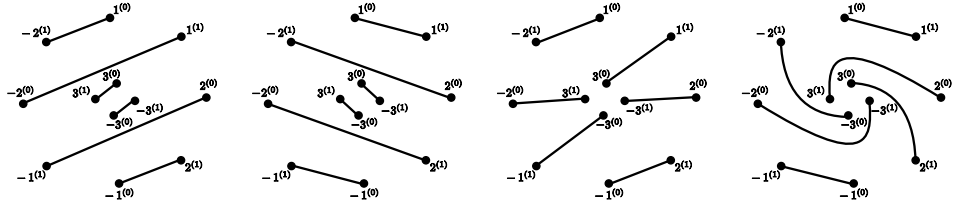
 
    \centering
    \begin{tabular}{cccc}
      \input{handshakeD3_1.tex} & \input{handshakeD3_2.tex} &
      \input{handshakeD3_3.tex} & \input{handshakeD3_4.tex}
    \end{tabular}
    \caption{Four different noncrossing handshake configurations in $\T_{D_3}$.}
\label{fig:handshakeD}
  \end{figure}

  Define a matching of
  $$\{\pm 1^{(0)},\pm 1^{(1)},\ldots,\pm n^{(0)},\pm n^{(1)}\}$$
  to be noncrossing of type $D_n$ if the points $\{\pm 1^{(0)},\pm 1^{(1)},\ldots,\pm (n-1)^{(0)},\pm (n-1)^{(1)}\}$ are arranged clockwise on a circle as in type $C_{n-1}$ and the points $\{\pm n^{(0)},\pm n^{(1)}\}$ form a small counterclockwise oriented square in the center of the circle, and the matching does not cross in this sense. A \Dfn{noncrossing handshake configuration} $T$ \Dfn{of type} $D_n$ is a noncrossing matching $T$ of type $D_n$, with the additional properties that $(i^{(1)},j^{(0)}) \in T$ if and only if $(-i^{(1)},-j^{(0)}) \in T$ and that the size of
  $$M_\pm := \{(i^{(1)},j^{(0)}) \in T : \mbox{$i$ and $j$ have opposite signs}\}$$
  is divisible by $4$. See Figure~\ref{fig:handshakeD} for examples of noncrossing handshake configurations of type $D_3$.

  As in the other types, we keep in mind that we could replace the linear Coxeter element by any Coxeter element to obtain labellings for the vertices of a noncrossing handshake configuration of type $D$.

  Define the Kreweras complementation $\K$ on $D_n$ by rotating the labels of the outer circle counterclockwise and the labels of the inner circle clockwise; more precisely, let $\kappa(i^{(0)}) := i^{(1)}$ and 
  \begin{eqnarray}
    \kappa(i^{(1)}) :=
      \left\{ \begin{array}{cl}
        (i+1)^{(0)}  & \mbox{if $i \in [n-2]$} \\
        (i-1)^{(0)}  & \mbox{if $i \in [-n+2]$} \\
        (-1)^{(0)}  & \mbox{if $i = n-1$} \\
        1^{(0)}      & \mbox{if $i = -n+1$} \\
        (-n)^{(0)}  & \mbox{if $i = n$} \\
        n^{(0)}      & \mbox{if $i = -n$.}
      \end{array} \right. \label{eq:krewerasD}
    \end{eqnarray}
  Then $(i^{(1)},j^{(0)}) \in T$ if and only if $\big(\kappa(j^{(0)}),\kappa(i^{(1)})\big) \in \K(T)$. To see this, observe that the only outer vertices changing sign are $\pm (n-1)^{(1)}$, and the only two inner vertices are $\pm n^{(1)}$. Thus, the size of $M_\pm$ for $\K(T)$ is again divisible by $4$. As an immediate consequence of the construction in \cite{AR2004}, we obtain that the map $\phi_{D_n} : \T_{D_n} \tilde{\longrightarrow} \NC(D_n)$ defined in the same way as for $\NC_n$ is well-defined and a bijection between noncrossing handshake configurations of type $D_n$ and $\NC(D_n)$.
  \begin{proposition} \label{prop:KrewerasD}
    The bijection $\phi_{D_n} : \T_{D_n} \tilde{\longrightarrow} \NC(D_n)$ is compatible with the Kreweras complementation, i.e., for $T \in \T_{D_n}$,
    $$\phi_{D_n}(\K(T)) = \K(\phi_{D_n}(T)).$$
  \end{proposition}
  \begin{proof}
    Let $\big(i^{(1)},j^{(0)}\big) \in T$. This implies that $\big(\kappa(j^{(0)},\kappa(i^{(1)}))\big) \in \K(T)$. Therefore, by checking the different cases in (\ref{eq:krewerasD}), we obtain $\phi_{D_n}(\K(T))\phi_{D_n}(T) = c$,
    and moreover, $\phi_{D_n}(\K(T)) = c \phi_{D_n}(T)^{-1} = \K(\phi_{D_n}(T))$.
  \end{proof}
  \begin{proposition} \label{prop:CSPNCD}
    The triple $\big(\NC(D_n),\Cat(D_n;q),\langle \K \rangle \big)$ exhibits the CSP.
  \end{proposition}
  \begin{proof}
    The $q$-Catalan number $\Cat(D_n;q)$ is given by
    \begin{align*}
      \Cat(D_n,q) = \begin{bmatrix} 2n-1 \\ n \end{bmatrix}_{q^2} + q^n \begin{bmatrix} 2n-2\\ n \end{bmatrix}_{q^2}.
    \end{align*}
    Let $d$ be an integer such that $d \big| 4(n-1)$ and let $\zeta$ be a primitive $d$-th root of unity. Then it follows again from \cite[Lemma 3.2]{EF2008} that $\Cat(D_n,q)$ reduces for $q = \zeta$ to
    \begin{eqnarray*}
      \big[\Cat(D_n,q)\big]_{q=\zeta} = 
        \left\{ \begin{array}{cl} 
          {\displaystyle\Cat(D_n)} &\mbox{if $d=1$} \\[6pt]
          {\displaystyle\Cat(D_n)} &\mbox{if $d=2$, $n$ even} \\[6pt]
          {\displaystyle \Cat(C_{n-1})} &\mbox{if $d=2$, $n$ odd} \\[6pt]
          {\displaystyle \Cat(C_{n/2})} &\mbox{if $d=4$, $4\big|n$} \\[6pt]
          {\displaystyle \Cat(C_{2(n-1)/d})} &\mbox{if $d \geq 4$ even, $d\big|2(n-1)$} \\[6pt]
          {\displaystyle \Cat(C_{(n-1)/d})} &\mbox{if $d \geq 3$ odd} \\[6pt]
          {\displaystyle 0} & \mbox{otherwise}
        \end{array} \right. \label{eq:CatD}
    \end{eqnarray*}
    \begin{figure} 
      \centering
      \begin{tabular}{c}
        \input{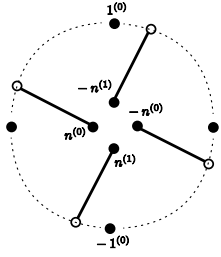}
      \end{tabular}
      \caption{A typical situation in $\T_{D_n}$ with $4$-fold symmetry for $4 = d \big|n$.} \label{fig:handshakeDsym}
    \end{figure}
    
    For $d=1$, this is obvious.
    
    For $d=2$, $n$ even, the symmetry property implies that $\K^{2(n-1)}(T) = T$ for all $T \in \T_{D_n}$.
    
    For $d=2$, $n$ odd, observe that $T \in \T_{D_n}$ is invariant under $2$-fold symmetry, i.e., $\K^{2(n-1)}(T) = T$ if and only if $\{\pm n^{(0)},\pm n^{(1)}\}$ forms a sub-matching of $T$. Therefore, those are counted by $\Cat(C_{n-1})$.
    
    For $d=4\big|n$, we want that $\K^{n-1}(T) = T$ and therefore, $\{\pm n^{(0)},\pm n^{(1)}\}$ must not form a sub-matching of $T$ and we are in a situation as indicated in Figure~\ref{fig:handshakeDsym}. This gives
    \begin{align*}
      \big|\{T \in \T_{D_n} : \K^{n-1}(T) = T\}\big| &= 2(n-1)\Cat(A_{(n-2)/2}) \\
                                                     &= \frac{4(n-1)}{n}\binom{n-2}{(n-2)/2} = \binom{n}{n/2},
    \end{align*}
    where the first $2$ comes from the $2$-fold rotation of the inner square, the $n-1$ is the number of possible connections between the inner square and the circle, and $\Cat(A_{(n-2)/2})$ is the number of noncrossing handshake configurations of the $n-2$ free points on the outer circle.
    
    For $d \geq 4$ even, $d\big|2(n-1)$, we have again that $\{\pm n^{(0)},\pm n^{(1)}\}$ forms a sub-matching of $T$ and we have immediately that
    $$\big|\{T \in \T_{D_n} : \K^{4(n-1)/d}(T) = T\}\big| = \Cat(C_{2(n-1)/d}).$$
    
    For $d \geq 3$ odd, it follows that $d\big|n-1$ and the same argument as in the previous case applies.
    
    The only otherwise case which is left is the case $d \geq 4$ even, $d \nmid 2(n-1)$. In this case, we see that $4 \big| d$ and it follows together with the symmetry property that there does not exist a $T \in \T_{D_n}$ such that $\K^{4(n-1)/d}(T) = T$.
  \end{proof}

\subsection{Type $I_2(k)$}

  For the dihedral groups, we obtain the theorem by straightforward computations. Let $I_2(k) = \langle a,b \rangle$ for two given simple reflections $a,b$ and fix the linear Coxeter element $c := ab$. Then $\NC(I_2(k))$ contains $\id,c$ and all $k$ reflections contained in $I_2(k)$. 

  \begin{proposition}
    The triple $\big(\NC(I_2(k)),\Cat(I_2(k);q),\langle \K \rangle \big)$ exhibits the CSP.
  \end{proposition}
  \begin{proof}
    The Kreweras complementation $\K$ on $\NC(I_2(k))$ has $2$ orbits, one is $\{\id,c\}$ and the other contains all $k$ reflections. On the other hand,
    \begin{align*}
      \Cat(I_2(k);q) &= \frac{[k+2]_q[2k]_q}{[2]_q[k]_q} \\
                     &= 
      \left\{ \begin{array}{cl} 
        (1+q^2+\dots+q^k)(1+q^k)                    &\mbox{if $k$ even} \\[6pt]
        1+q^2+\dots+q^{k-1}+q^k+q^{k+1}+\dots+q^{2k}  &\mbox{if $k$ odd},
      \end{array} \right.
    \end{align*}
    and the proposition follows.
  \end{proof}

\subsection{Exceptional types} \label{sec:exceptionaltypes}

    For the exceptional Coxeter groups,
    $$\Cat(W;q) \mod(q^{2h}-1)$$
    can be simply computed and by (\ref{eq:CSPorbitlengths}), we need to find the following orbit lengths, where $i*j$ is shorthand for $i$ orbits of length $j$:
    \begin{align*}
      F_4 &: 8*12, 1*4, 1*3, 1*2, \\
      H_3 &: 3*10, 1*2, \\
      H_4 &: 9*30, 1*5, 1*3, 1*2, \\
      E_6 &: 30*24, 8*12, 1*8, 1*4, 1*3, 1*2, \\
      E_7 &: 230*18, 3*6, 1*2, \\
      E_8 &: 832*30, 5*15, 3*10, 2*5, 1*3, 1*2.
    \end{align*}
    Those orbit lengths were verified with a computer; as mentioned above, they can be deduced as well from \cite{KM2008}.
    
\section{The Panyushev CSP for nonnesting partitions} \label{sec:CSPNN}

\begin{figure}
    \centering
    \begin{tabular}{c}
      \input{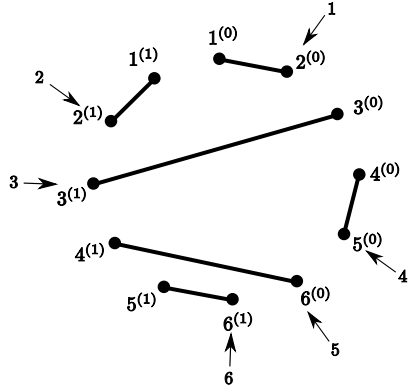}
    \end{tabular}
    \caption{The nonnesting labels on a noncrossing handshake configuration in $\T_6$.}
    \label{fig:handshakeNNA}
  \end{figure}

   In this section, we prove Theorem~\ref{th:CSPNC}(ii) for every type individually by providing a bijection between nonnesting partitions and noncrossing handshake configurations which maps the Panyushev complementation to the Kreweras complementation. We consider the same noncrossing handshake configurations as before,
but we use a different labelling to refer to the vertices.  In type $\An$, 
we label the vertices on the outer circle by $\{1^{(0)},\ldots,n^{(0)},n^{(1)},\ldots,1^{(1)}\}$ in clockwise order.  E.g., the noncrossing handshake configuration shown in Figure~\ref{fig:handshakeA} is relabeled as shown in Figure~\ref{fig:handshakeNNA}.
  
\subsection{Type $A$}
    
\begin{figure}
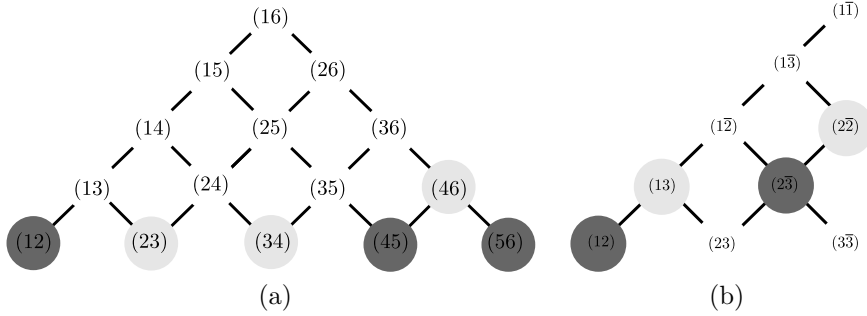

    \centering
    \begin{tabular}{cc}
      \input{rootposetA5.tex}
      &
      \input{rootposetB3.tex} \\
      (a) & (b)
    \end{tabular}
    \caption{(a) An antichain and its image under the Panyushev complementation in the root poset of type $A_5$; (b) another antichain and its image in the root poset of type $C_3$.}
    \label{fig:rootposets}
  \end{figure}

  Let $\Phi^+ := \{ (i,j) = e_i - e_j : 1 \leq i < j \leq n\}$ be the set of all transpositions identified with a set of positive roots for $\An$. The root poset structure on $\Phi^+$ is given by
  \begin{align}
    (i,j) \leq (i',j') \Leftrightarrow i' \leq i < j \leq j', \label{eq:rootrelationA]}
  \end{align}
  see Figure~\ref{fig:rootposets}(a) for an example.
  Let $I = \{(i_1,j_1),\ldots,(i_k, j_k)\} \in \NN(\An)$ such that $i_1 < \dots < i_k$. Observe that (\ref{eq:rootrelationA]}) implies $j_1 < \dots < j_k$ as well. Define a map
  $$\psi_{\An} : \NN(\An) \longrightarrow \T_n$$
  as follows: for $1 \leq \ell \leq k$, mark the vertex $j_\ell^{(0)}$ with $i_\ell$ and for $i \in [n] \setminus \{i_1,\ldots,i_k\}$ mark the vertex $i^{(1)}$ with $i$. Now, for $1 \leq i \leq n$, in increasing order, match the vertex marked with $i$ with the first non-matched vertex, where first is interpreted counterclockwise from the marked vertex if $i \in \{i_1,\ldots,i_k\}$ and clockwise 
from the marked vertex if $i \notin \{i_1,\ldots,i_k\}$. For example, for the antichain
  $$I = \big\{(1,2),(4,5),(5,6) \big\} \in \NN(A_5)$$
  considered in Figure~\ref{fig:rootposets}(a), we have $\psi_{\An}(I) = T$, where $T \in \T_6$ is the noncrossing handshake configuration shown in Figures~\ref{fig:handshakeA} and \ref{fig:handshakeNNA}.
  
  To show that $\psi_{\An}$ is a bijection, we now define its inverse map $\psi'_{\An}: \T_n \longrightarrow \NN(\An)$. Let $T \in \T_n$. Mark all $j^{(\beta)}$ for which $(i^{(\alpha)},j^{(\beta)}) \in T$ with $i < j$, or with $i = j$ and $(\alpha,\beta) = (0,1)$. Next, label all marks $i^{(1)}$ with $i$, and then label all marks $i^{(0)}$ clockwise with the remaining labels in $[n]$. The antichain $\psi'_{\An}(T)$ is then given by
  $$\psi'_{\An}(T) = \big\{(i,j) : \text{ vertex } j^{(0)} \text{ is marked by } i \big\}.$$
  
  \begin{proposition}
    The map $\psi'_{\An}$ is well-defined and the inverse of $\psi_{\An}$. In particular, $\psi_{\An}: \NN(\An) \tilde{\longrightarrow} \T_n$ is a bijection.
  \end{proposition}
  \begin{proof}
    To see that $\psi'_{\An}$ is well-defined, we have to check that any marked vertex $j^{(0)}$ is marked with some $i<j$. Assume that $j^{(0)}$ is marked with $j$. This implies that the set $\{1^{(0)},\ldots,(j-1)^{(0)},(j-1)^{(1)},\ldots,1^{(1)}\}$ contains $j-1$ marked vertices and forms therefore a sub-matching -- a contradiction to the fact that $j$, as it is marked, is matched to some element in this set.
    
    As in the process of applying $\psi'_{\An}$ and of applying $\psi_{\An}$ the same vertices get marked, $\psi'_{\An}$ is in fact the inverse of $\psi_{\An}$.
  \end{proof}
  \begin{theorem} \label{th:PanyushevKreweras}
    The bijection $\psi_{\An}$ is compatible with the Panyushev respectively the Kreweras complementation. For $I \in \NN(\An)$, we have
    $$\K(\psi_{\An}(I)) = \psi_{\An}(\P(I)).$$
  \end{theorem}
  To prove this theorem, we first have to understand how the Panyushev complementation behaves in type $A$. Recall that the support $\supp(I)$ of some antichain $I \in \NN(\An)$ is given by $\supp(I) := \bigcup_{(i,j) \in I} \{s_i,\ldots,s_{j-1}\}$. Next, set
  $$\hat{I} = \big\{(i'_1,j'_1),\ldots,(i'_k,j'_k)\}:= I \cup \{(i,i) : s_{i-1},s_i \notin \supp(I) \big\}$$
  such that $i'_1 < \ldots < i'_k$, where the dummies $s_0,s_n$ are supposed not to be in $\supp(I)$. The Panyushev complementation is then given by
  $$\P(I) = \big\{(i'_2-1,j'_1+1),\ldots,(i'_k-1,j'_{k-1}+1)\big\} \in \NN(\An).$$
  \begin{proposition} \label{prop:supporthandshake}
    Let $I$ be a nonnesting partition. Then $s_k \notin \supp(I)$ if and only if $\{i^{(0)},i^{(1)} : 1 \leq i \leq k \}$ defines a submatching of $\psi_{\An}(I)$. In particular,
    $$(i^{(0)},i^{(1)}) \in \psi_{\An}(I) \Leftrightarrow (i,i) \in \hat{I}.$$
  \end{proposition}
  \begin{proof}
    The proposition follows directly from the definition.
  \end{proof}
  \begin{example}\label{ex:NNsupport}
    The noncrossing handshake configuration $T$ in Figure~\ref{fig:handshakeNNA} is the image of $I = \{(1,2),(4,5),(5,6)\} \in \NN(A_5)$ under $\psi_{A_5}$. The complement of the support of $I$ is $S \setminus \supp(I) = \{s_2,s_3\}$.  The submatchings guaranteed by the Proposition are those of the form $\{1^{(0)},1^{(1)},\ldots,k^{(0)},k^{(1)}\}$ for $k \in \{ 2, 3\}$.
  \end{example}
  \begin{proof}[Proof of Theorem~\ref{th:PanyushevKreweras}]
    As it is easier to see, we describe the analogous statement for $\psi'_{\An}$. $\psi'_{\An}(\K(T))$ can be described in terms of $\psi'_{\An}(T)$ as follows: a marked $i^{(0)}$ is turned to a marked $(i+1)^{(0)}$ (unless $i=n$ when the mark disappears), and for a marked $i^{(1)}$, we obtain a marked $(i-1)^{(1)}$ (unless $i=1$ when the mark disappears). If $(i^{(0)},i^{(1)}) \in T$, the marked $i^{(1)}$ is replaced by a marked $(i+1)^{(0)}$. The theorem follows with Proposition~\ref{prop:supporthandshake} and the description of $\P(I)$ in terms of $\hat{I}$.
  \end{proof}

\subsection{Types $B$ and $C$}

  In contrast to the situation for reflection groups, the notion of the root system does not coincide for types $B$ and $C$. The resulting root posets turn out to be isomorphic (as posets) but not equal. 
Thus, it suffices to study the Panyushev complementation on one of the two. 
As the connection between the root poset of type $C_n$ and the root poset of type $A_{2n-1}$ is straightforward, whereas there is a little more work to do in type $B_n$, we will study nonnesting partitions of type $C_n$. This corresponds to the fact that the type $C_n$ Dynkin diagram can be obtained from the type $\An$ Dynkin diagram through a \lq\lq folding process\rq\rq.

  The set of reflections identified with a set of positive roots in type $C_n$ is given by
  $$\Phi^+ := \{ (i,j) = e_i - e_j : 1 \leq i < j \leq n \} \cup \{(i,\overline{j}) = e_i + e_j : 1 \leq i \leq j \leq n\}.$$
  See Figure~\ref{fig:rootposets}(b) for the root poset of type $C_3$ as an example.
  
  To understand nonnesting partitions of type $C_n$, observe that an antichain in $\Phi^+$ can be identified with a symmetric antichain in the root poset of type $A_{2n-1}$: there is an involution $\delta$ on $\NN(\An)$ by horizontally flipping the root poset of type $\An$, i.e., replacing the positive root $(i,j)$ by $(n+1-j,n+1-i)$. In other words, $\delta$ is the induced map coming from the involution on the Dynkin diagram sending one linear ordering to the other. Define an antichain $I \in \NN(\An)$ to be \Dfn{symmetric} if it is invariant under this involution.
  It is well-known that $\NN(C_n)$ can be seen as the set of all antichains $A \in \NN(A_{2n-1})$ which are symmetric,
  $$\NN(C_n) \cong \big\{I \in \NN(A_{2n-1}) : \delta(I) = I \big\}.$$
  Moreover, this identification is compatible with the Panyushev complementation,
  $$\delta(I) = I \Leftrightarrow \delta(\P(I)) = \P(I).$$
  This allows us to study this complementation on nonnesting partitions of type $C_n$ in terms of symmetric nonnesting partitions of type $A_{2n-1}$.
  
  On the other hand, we have seen above that the bijection $\phi_{A_{2n-1}} : \T_{2n} \longrightarrow \NC(A_{2n-1})$ restricts to a bijection $\phi_{C_n} : \T_{C_n} \longrightarrow \NC(C_n)$. Therefore, we want to show that the bijection $\psi_{A_{2n-1}} : \NN(A_{2n-1}) \longrightarrow \T_{2n}$ gives rise to a bijection $\psi_{C_n} : \NN(C_n) \longrightarrow \T_{C_n}$ which is again compatible with the Panyushev and the Kreweras complementation.
  \begin{lemma} \label{lem:Binvolution}
    The involution $\delta$ on $I$ for $I \in \NN(A_{n-1})$ can be described in terms of the Kreweras complementation as
    $$\psi_{A_{n-1}}(\delta(I)) = \K^{n}(\psi_{A_{n-1}}(I)).$$
  \end{lemma}
  \begin{proof}
     For $T \in \T_n$, we have
    $$(i^{(\alpha)},j^{(\beta)}) \in T \Leftrightarrow \big((n+1-j)^{(\beta^c)},(n+1-i)^{(\alpha^c)}\big) \in \K^n(T),$$
    where $\alpha,\beta \in \{0,1\}$ and $\alpha^c$ (resp. $\beta^c$) denotes the complement of $\alpha$ (resp. $\beta$) in $\{0,1\}$. It is straightforward to check that this observation implies that
    $$\psi'_{A_{n-1}}\big(\K^{n}(\psi_{A_{n-1}}(I))\big) = \delta(I).$$
  \end{proof}
  \begin{theorem}
     $\psi_{A_{2n-1}}$ restricts to a well-defined bijection $\psi_{C_n} : \NN(C_n) \longrightarrow \T_{C_n}$.
  \end{theorem}
  \begin{proof}
    The statement of the theorem is equivalent to the statement that
    $$\delta(I) = I \Leftrightarrow \K^n(\psi_{A_{2n-1}}(I)) = \psi_{A_{2n-1}}(I).$$
    This follows directly from the previous lemma.
  \end{proof}

\subsection{Type $D$}

Fix the numbering of the Dynkin diagram of type $D_n$ so that $n-2$ is adjacent to $n-1$, $n$, and $n-3$.
We consider the involution
$\delta$
of this diagram which interchanges $n$ and $n-1$.  It acts on $\NN(D_n)$, $\NC(D_n)$, and $\T_{D_n}$.  On $\T_{D_n}$, it acts by rotating the inner four 
vertices by a half turn.  It is convenient to define a new type of 
noncrossing handshake configuration, which we denote $\T_{D_n/\delta}$: this consists of 
$4n-4$ external vertices, labelled as in a $C_{n-1}$ noncrossing handshake configuration, such
that either all the vertices participate in a $180^\circ$-rotationally symmetric
noncrossing matching 
(in which case we simply have a $C_{n-1}$ noncrossing handshake configuration) or else all
but four vertices participate in a $180^\circ$-rotationally symmetric
noncrossing matching, while the four remaining vertices are isolated but
have the property that any two of them could be attached without creating
any crossings.  
It is clear that elements of $\T_{D_n/\delta}$ correspond to $\delta$-orbits in $\T_{D_n}$.

\subsubsection{Defining a map from $\NN(D_n)/\delta$ to $\mathcal T_{D_n/\delta}$}
Note that $\rot$ acts naturally on $\T_{D_n/\delta}$, while $\Pan$ acts naturally
on $\delta$-orbits in $\NN(D_n)$.  We will begin by showing that 
$( \T_{D_n/\delta},\rot )$ and $( \NN(D_n)/\delta,\Pan )$ are isomorphic as sets
with a cyclic action.  

In this subsection, we will define a cardinality-preserving 
bijection from $\delta$-orbits in
$\NN(D_n)$ to $\T_{D_n/\delta}$ which we will denote by
$\psi_{D_n/\delta}$.  (In fact, for notational convenience, we will write 
$\psi_{D_n/\delta}$ as a map from $\NN(D_n)$ to $\T_{D_{n}/\delta}$ 
which is constant on $\delta$-orbits.)    
We will then show that it is possible to refine $\psi_{D_n/\delta}$  to a 
bijection from $\NN(D_n)$ to $\T_{D_n}$.  

{\em Singleton $\delta$-orbits in $\T_{D_n/\delta}$.}  Such an element consists
of a type $C_{n-1}$ noncrossing handshake configuration on $4n-4$ external vertices 
$1^{(0)},\dots,(2n-2)^{(0)},(2n-2)^{(1)},\dots
1^{(1)}$.  

{\em Singleton $\delta$-orbits in $\NN(D_n)$.}  Such an element of 
$\NN(D_n)$ corresponds to a single element of $\NN(B_{n-1})$.  We 
reinterpret this as an element of $\NN(C_{n-1})$,  which corresponds
(as we have already seen) to an element of $\NN(A_{2n-3})$ fixed under
the involution of the $A_{2n-3}$ diagram.

{\em Map from singleton $\delta$-orbits in $\NN(D_n)$ to $\T_{D_n}$.} 
We define $\psi_{D_n/\delta}$ on a singleton $\delta$-orbit by sending
the type 
$A_{2n-3}$ antichain to an $A_{2n-3}$ noncrossing handshake configuration, using
$\psi_{A_{2n-3}}$.  

\medskip

Now we consider the doubleton $\delta$-orbits.  
Write $H$ for the $2n-2$ vertices 
$\{(n-2)^{(1)},\dots,1^{(1)},1^{(0)},\dots,n^{(0)}\}$, and $H^c$ for the
other $2n-2$ vertices on the boundary.  

{\em Doubleton $\delta$-orbits in $\T_{D_n}$.}  These correspond to 
elements of $\T_{D_n/\delta}$ which have four vertices of degree zero.

{\em Doubleton $\delta$-orbits in $\NN(D_n)$.}  Let $I$ be an antichain
in such an orbit.  Write $\overline I$ for the collection of type
$A_{2n-3}$ roots obtained by taking each root in $I$, passing first to 
$B_{n-1}$, identifying the root poset of $B_{n-1}$ with that of $C_{n-1}$, 
and then unfolding to one or two roots in $A_{2n-3}$.  
Note that $\overline I$ is typically not an antichain.

\begin{example}
Consider the $D_n$ antichain consisting of $\alpha_n+\alpha_{n-2}$
and $\alpha_{n-1}$.  The former contributes elements $(n-1,n+1)$ and 
$(n-2,n)$, while the latter contributes $(n-1,n)$.  This does not form
an antichain.  There will often be two elements in $\overline I$ with 
first co-ordinate $n-1$, and two elements with second co-ordinate $n$.
\end{example}

We also associate to $I$ an antichain in $\Phi_{A_{2n-3}}$, defined as follows.  
Consider the elements of 
$\overline I$ which lie in the square with opposite corners at $(1,2n-2)$ and 
$(n-1,n)$. (We call this square $R$.)
Record the first coordinates of these as $i_1,\dots,i_r$,
and the last as $j_1,\dots,j_r$.  

Note that $j_1=j_2$ and $i_r=i_{r-1}$ are
possible (occurring when $\overline I$ is not an antichain).  
Define $\widehat I$ by
replacing these $r$ elements of $\overline I$ by the $r-1$ elements
$(i_1,j_2), (i_2,j_3),\dots, (i_{r-1},j_r)$.  (In the case that $r=1$, 
the result is that $\widehat I\cap R=\emptyset$.)

{\em The map from doubleton $\delta$-orbits in $\NN(D_n)$ to 
doubleton $\delta$-orbits in $\NC(D_n)$.}  We define 
$\psi_{D_n/\delta}(I)$ in several steps.  Using Lemma \ref{widehatisgood},
below, we know that $\widehat I\in \NN(C_{n-1})$.  Therefore, we can
consider $\psi_{C_{n-1}}(\widehat I)\in \mathcal T_{C_{n-1}}$.  
Lemma \ref{lem3} below guarantees that there
are at least two edges in this diagram which run from
vertices in $H$ to vertices in $H^c$.  Remove the two such edges which are closest
to the center.  The result is a noncrossing handshake configuration of type $D_n/\delta$
as defined above.  This is $\psi_{D_n/\delta}(I)$.   

\subsubsection{Defining $\psi_{D_n}$}
We now consider refining $\psi_{D_n/\delta}$ to a map from $\NN(D_n)$ to
$\mathcal T_{D_n}$.  

We use the convention that a type 
$D$ noncrossing handshake configuration has the same outside labels as for 
type $D/\delta$ noncrossing handshake configurations, with four internal vertices which 
are numbered by congruence classes modulo 4,
increasing in counter-clockwise order.    
We count as ``positive'', external vertices with label $(0)$, and 
the internal vertices $0$ and $3$, and as ``negative'', external vertices
with the label $(1)$ and the internal vertices $1$ and $2$.  
In a noncrossing handshake configuration, the number of edges that connect a positive
vertex to a negative vertex must be divisible by 4.  

If a noncrossing handshake configuration $T$ of type $D_n/\delta$ has no isolated vertices,
this requirement means that there is a unique way of completing $T$ 
to a type $D_n$ configuration, while if $T$ has four isolated vertices,
then there are two ways of completing $T$ to a type $D_n$ configuration.  

For $a,b$ outer vertices, write $d(a,b)$ for the clockwise distance from 
$a$ to $b$.  Write $e_I(a,b)$ for the number of vertices in the clockwise
interval from $a$ to $b$, including $b$ but not $a$, and which are not
on the clockwise end of an edge in $\psi_{D_n/\delta}(I)$.  

For $I$ an antichain in $\NN(D_n)$ in a doubleton $\delta$-orbit, 
define $s(I)$ to be $0$ if the 
root of $I$ whose image in $\overline I$ is $(i,n)$ with $i$ as small
as possible, has $\alpha_{n-1}$ in its support; otherwise, set $s(I)=1$.  

We now define $\psi_{D_n}(I)$.  If $I$ is in a singleton $\delta$-orbit,
then define $\psi_{D_n}(I)$ to be $\psi_{D_n/\delta}(I)$ together with
edges connecting the internal vertices in the unique possible way.  

If $I$ is in a doubleton $\delta$-orbit, define $\psi_{D_n}(I)$ by starting
with $\psi_{D_n/\delta}(I)$ and, for each singleton external vertex $v$, 
attach it to the internal vertex whose number
is given by:
$n-d(v,(n-1)^{(0)})+2s(I)+2e_I(v,(n-1)^{(0)})$.  

\begin{example}
	For the root poset of type $D_3$ with simple roots
	$$\alpha_1 = e_1 - e_2, \alpha_2 = e_2 - e_3, \alpha_3 = e_2 + e_3,$$
	the four antichains $\emptyset, \{\alpha_1,\alpha_2,\alpha_3\}, \{\alpha_2\}, \{\alpha_1,\alpha_3\}$ are mapped by $\psi_{D_n}$ to the four noncrossing handshake configurations in $\T_{D_3}$ shown in Figure~\ref{fig:handshakeD} from left to right.
\end{example}

\subsubsection{Proof that $\psi_{D_n}$ is well-defined and is a bijection}
There are several lemmas
which must be  established to show that the definition given above makes sense, and yields a bijection.  

\begin{lemma} \label{widehatisgood}
$\widehat I$ is in $\NN(A_{2n-3})$.  Further, the map from 
$I$ to $\widehat I$ is injective, and its image consists of all the antichains 
in $\NN(C_{n-1})$ (thought of as a subset of $\NN(A_{2n-3})$)
except those containing $(n-1,n)$.  
\end{lemma}

\begin{proof}
The inverse map is clear, since $i_r$ must be $n$ and $j_1$ must be
$n-1$.  This inverse map can be applied to any antichain in
$\NN(C_{n-1})$ except those containing $(n-1,n)$.  
\end{proof}

Now, since $\widehat I$ is in $\NN(C_{n-1})$, its image under the
bijection $\psi_{A_{2n-3}}$ is a  type $C_{n-1}$ noncrossing handshake configuration.  
The following lemma is useful.  

\begin{lemma}\label{lem2}
The image of $\psi_{C_{n-1}}$ applied to antichains with no roots in $R$,
consists exactly of those type $C_{n-1}$ noncrossing handshake configurations with no edges from 
$\{(n-1)^{(1)},\dots,1^{(1)}, 1^{(0)},\dots,(n-1)^{(0)}\}$ to the other vertices.  
\end{lemma}

\begin{proof}
The first $n-1$ edges in the noncrossing handshake configuration will all connect
vertices in $\{(n-1)^{(1)},\dots,1^{(1)}$, $1^{(0)},\dots,(n-1)^{(0)}\}$, 
which uses up
all those vertices.  
\end{proof}

\begin{lemma} \label{lem3}
The image of $\psi_{C_{n-1}}$ applied to $\widehat I$ for $I \in \NN(D_n)$, 
consists of exactly
those type $C_{n-1}$ 
noncrossing handshake configurations with the property that there is at least 
one edge (and therefore at least two edges) from $H$ to $H^c$.  
\end{lemma}

\begin{proof}
We have already shown that as $I$ runs through 
$\NN(D_n)$, we have that $\widehat I$ runs through those antichains
in $\NN(C_{n-1})$ not containing $(n-1,n)$.
The image under $\Pan^{-1}$ of type $C_{n-1}$ antichains not containing 
$(n-1,n)$ is exactly the $C_{n-1}$ antichains whose intersection with
$R$ is non-empty.  
Now apply Lemma \ref{lem2} to $\Pan^{-1}(\widehat I)$, together
with the fact that $\rot \circ \psi_{C_{n-1}} = \psi_{C_{n-1}}\circ \Pan$.  
\end{proof} 

We now have the pieces in place to establish the following proposition:

\begin{proposition} The map $\psi_{D_n/\delta}$ is a bijection from
$\NN(D_n/\delta)$ to $\T_{D_n/\delta}$.  
\end{proposition}

\begin{proof} It is clear that $\psi_{D_n/\delta}$ takes singleton $\delta$
orbits in $\NN(D_n)$ bijectively to the noncrossing handshake configurations in 
$\mathcal T_{D_n/\delta}$ which contain no isolated vertices.  It is also
clear that $\psi_{D_n/\delta}$ is an injection from doubleton orbits in
$\NN(D_n)$ into the $\mathcal T_{D_n/\delta}$ noncrossing handshake configurations with
four isolated vertices.  Finally, given such a diagram, there is a 
unique way to reattach the isolated vertices to obtain a 
$\T_{C_{n-1}}$ noncrossing handshake configuration such that the reattached edges cross
from $H$ to $H^c$.  It follows that $\psi_{D_n/\delta}$ is a bijection.
\end{proof}

We now proceed to show that $\psi_{D_n}$, as defined above, is a bijection
from $\NN(D_n)$ to $\T_{D_n}$.  To begin with, we need the following 
lemma which gives a condition equivalent to the parity condition on
the number of edges in a type $D_n$ noncrossing handshake configuration which connect
positive and negative vertices.  

\begin{lemma} \label{restatement}
The condition that the number of edges joining a positive vertex
to a negative vertex be divisible by four, is equivalent to the condition that
a positive, even-numbered singleton vertex must
be connected to an internal vertex of odd parity, and similarly for 
the other possible choices of singleton vertex, where changing either
``positive'' or ``even-numbered'' reverses the parity of the internal
vertex.
\end{lemma}

We are now ready to prove that $\psi_{D_n}$ is a bijection.  

\begin{lemma} $\psi_{D_n}$ is a bijection from $\NN(D_n)$ to 
$\mathcal T_{D_n}$.  
\end{lemma}
\begin{proof} We must show that if $v$ and $v'$ are singleton
vertices in $\psi_{D_n/\delta}(I)$, such that the next singleton vertex after
$v$ in counter-clockwise order is $v'$, then 
the vertex to which $v'$ is
attached is one step counter-clockwise from that to which $v$ is attached.  
We evaluate $-d(v',(n-1)^{(0)})+d(v,(n-1)^{(0)})=-d(v',v)$ by
counting the vertices between $v'$ and $v$ (including $v$ but not
$v'$).  
Each edge on the outer rim between $v$ and $v'$ contributes
$-2$ to $-d(v',v)$ (one for each of its endpoints), 
and also contributes 2 to $2e_I(v',(n-1)^{(0)})-2e_I(v,(n-1)^{(0)})=2e_I(v',v)$.
The only other contribution to $2e_I(v',v)$ is an 
additional 2 coming from the vertex $v$, and also 
$-d(v',v)$ has an additional $-1$ coming from $v$.  Thus
the total effect is that $v'$ is attached one step counter-clockwise from 
$i$.  

The condition provided by Lemma \ref{restatement} is also clear from
the definition.  (Note that the complicated terms don't have any effect
on the parity of the vertex to which we connect $v$.)

Bijectivity follows from bijectivity for $\psi_{D_n/\delta}$ together with
the fact that the two elements of a doubleton $\delta$ orbit in 
$\NN(D_n)$ will be mapped to different noncrossing handshake configurations.
\end{proof}

\subsubsection{Compatibility between Panyushev complementation and rotation}

We will first prove that $\psi_{D_n/\delta}$ expresses the compatibility 
between Panyushev complementation for $\NN(D_n)/\delta$ and rotation of
$D_n/\delta$ noncrossing handshake configurations, and then we will prove the similar 
result for $\psi_{D_n}$.  

\begin{proposition}\label{modsigmacomp}
For $I \in \NN(D_n)$, we have that 
$$\psi_{D_n/\delta}(\Pan(I))=\rot(\psi_{D_n/\delta}(I)).$$
\end{proposition}

\begin{proof}
We consider three cases separately.  The first case is the case that 
$I$ is in a singleton $\delta$-orbit, 
in which case the result follows immediately
from the analogous result for type $C_{n-1}$.  

The second case is when $\widehat I\cap R \ne \emptyset$.  

\begin{lemma} If $\widehat I\cap R \ne \emptyset$, then 
$\widehat{\Pan(I)}=\Pan(\widehat I)$ and 
$\psi_{C_{n-1}}(\Pan(\widehat I)) =\rot(\psi_{C_{n-1}}(\widehat I))=
\rot(\psi_{D_n/\delta}(I))$
\end{lemma}

\begin{proof} The fact that {$\widehat{\Pan(I)}=\Pan(\widehat I)$} in this 
case follows from the definitions.
The compatibility of $\Pan$ and $\K$ in type $C$ implies that $\psi_{C_{n-1}}(\Pan(\widehat I)) =\rot(\psi_{C_{n-1}}(\widehat I))$.

Finally, we wish to show that $\rot(\psi_{C_{n-1}}(\widehat I))=
\rot(\psi_{D_n/\delta}(I))$.
The result which has to be established is that the pair of innermost edges
in $\rot(\psi_{C_{n-1}}(\widehat I))$ is the rotation of the innermost edges of 
$\psi_{C_{n-1}}(\widehat I)$ .  This is true because, in order for the
innermost edges no longer to be innermost, they must no longer run between
the two sides of the diagram.  But this would then imply that there were
no edges between $H$ and $H^c$ in $\psi_{C_{n-1}}(\widehat{\Pan(I)})$,
contrary to Lemma \ref{lem3}.
\end{proof}

We now consider the case that $\widehat I\cap R= \emptyset$.  In this
case, in contrast to the previous one, the proof does not pass 
through the similar statement in type $C$.  

Let $\widehat X=\Pan(\widehat I)$.  It is immediate from
the definition of Panyushev complementation that $\widehat X\cap R=
(n-1,n)$.  
By Lemma \ref{lem3} it follows that 
$\psi_{C_{n-1}}(\widehat X)$ has no edges from $H$ to $H^c$.  

By the compatibility of $\Pan$ and $\K$ in type $C$, we have that 
$\psi_{C_{n-1}}(\widehat X)=\rot (\psi_{C_{n-1}}(\widehat I))$.  
The innermost edges of $\psi_{C_{n-1}}(\widehat I)$ connecting 
$H$ to $H^c$, after rotation, no longer connect $H$ to $H^c$.  Thus,
in $\psi_{C_{n-1}}(\widehat X)$, those edges connect
$(n+1)^{(0)}$ to some
$z'$ in $H^c$ and
$(n-2)^{(1)}$ to some (symmetrical) $z$ in $H$.  

\begin{lemma} $\psi_{C_{n-1}}(\widehat{\Pan(I)})$ can be obtained from $\psi_{C_{n-1}}
(\widehat X)$ by removing the edges connected to $(n+1)^{(0)}$ and 
$(n-2)^{(1)}$ and replacing them by the other possible pair of symmetrical
edges.\end{lemma}

\begin{proof}
$\overline I \cap R$ necessarily equals $(n-1,n)$.
Let $Y=\Pan(I)$.  There are two possibilities for $\overline Y\cap R$:
it equals either $\{(n-2,n),(n-1,n),(n-1,n+1)\}$ or 
$\{(n-1,n)\}$,
depending on whether or not $\overline I$ has any entries on the 
$(n-1)$-th row (or equivalently the $(n+1)$-th column).   
The corresponding values of $\widehat Y\cap R$ are $\{(n-2,n),(n-1,n+1)\}$ and 
$\emptyset$.  

Now consider applying $\psi_{C_{n-1}}$ to $\widehat X$ and $\widehat Y$.  
Suppose first that we are in the case that $\widehat Y\cap R=\emptyset$.  This
means that the $(n-1)$-th row is empty in $\overline I$, so in $\widehat I$, both
$R$ and the row below $R$ are empty.  We have seen already that the fact
that $R$ is empty means that there are no edges between vertices numbered
at most $n-1$ and those numbered at least $n$.  A similar argument
shows that the absence of roots in the $(n-1)$-th row implies that the 
vertices numbered at most $n-2$ are connected to other vertices in that
set.  It follows that $(n-1)^{(0)}$ and $(n-1)^{(1)}$ are connected in
$\psi_{A_{2n-3}}(\widehat I)$.    
By symmetry, $n^{(0)}$ and $n^{(1)}$ are also.  

In determining $\psi_{C_{n-1}}(\widehat X)$, $n^{(0)}$ gets the label $n-1$.
In determining $\psi_{C_{n-1}}(\widehat Y)$, the label $n-1$ goes to $(n-1)^{(1)}$, the
symmetrically opposite vertex.  We know that $\psi_{C_{n-1}}(\widehat Y)$ has no
edges connecting vertices $\leq n-1$ with those $\geq n$, so the 
result of adding the $n-1$-th edge is to complete the matchings among
the vertices $\leq n-1$.  It follows that when we evaluate $\psi_{C_{n-1}}(\widehat X)$
instead, vertex $n^{(0)}$ will necessarily be connected to the
same vertex as $(n-1)^{(1)}$ was in $\psi_{C_{n-1}}(\widehat Y)$.  This means that,
while $n^{(0)}$ and $(n-2)^{(1)}$ are connected in 
$\psi_{C_{n-1}}(\widehat X)$, we have that $n^{(0)}$ and $(n+1)^{(0)}$ are connected
in $\psi_{C_{n-1}}(\widehat Y)$, establishing the claim. 

Now consider the case that $\widehat Y\cap R=\{(n-2,n),(n-1,n+1)\}$.  
In determining $\psi_{C_{n-1}}(\widehat Y)$, we have $n^{(0)}$ receives label
$n-2$ and $(n+1)^{(0)}$ receives label $n-1$.  Since $\widehat X$ and $\widehat Y$
only differ inside $R$, we have that the $n-2$-th column is empty
in $\widehat X$, so $(n-2)^{(1)}$ receives the $n-2$ label; and we also
have that $n^{(0)}$ receives the label $n-1$.  

Let us write $b$ for the vertex joined in $n^{(0)}$ in $\widehat X$, 
and $a$ for the vertex joined in $(n-2)^{(1)}$ in $\widehat X$.  Note 
that in $\widehat X$, there are no edges between $H$ and $H^c$, so,
prior to the $(n-2)$-th edge being drawn, the four available vertices
in $H$ are $(n-2)^{(1)},a,b,n^{(0)}$ (in clockwise order).  

Now consider what happens when we evaluate $\psi_{C_{n-1}}(\widehat Y)$.  When 
adding  the $(n-2)$-th edge, we connect $n^{(0)}$ to the next available
vertex counter-clockwise from it, which is $b$.  Next, we connect to
$(n+1)^{(0)}$ the next available vertex counter-clockwise from it, which
is $a$.  

The result is that $n^{(0)}$ is attached to the same vertex in 
$\widehat X$ and $\widehat Y$, but the vertex attached to $(n+1)^{(0)}$ in
$\widehat Y$ is attached to $(n-2)^{(1)}$ in $\widehat X$.  This suffices to
establish the claim.  \end{proof}

The final case of the proposition now follows, 
because the only edges between $H$ and $H^c$ in
$\psi_{C_{n-1}}(\widehat {\Pan(I)})$ are the new edges identified above, 
whose four end-vertices
are the result of rotating clockwise the four degree zero 
vertices of $\psi_{D_n/\delta}(I)$.
\end{proof}

In order to show the compatibility between $\psi_{D_n}$ and Panyushev 
complementation, 
we must study the relationship between $s(I)$ and $s(\Pan(I))$.  
It is straightforward
to check that $s(I)$ and $s(\Pan(I))$ are the same iff $I$ contains a root
supported over vertex $n-2$ but neither $n-1$ nor $n$.  This is equivalent
to saying that $\widehat I$ includes some root $(j,n-1)$ (i.e., a root
on the row just below $R$).  
This can also be described in terms of $\psi_{D_n}(I)$, as in the lemma below.

\begin{lemma} For $I$ an $A_{2n-3}$-antichain, $I$ contains a root 
$(j,n-1)$ iff $\psi_{D_n}(I)$ contains an edge joining $n-1$ to $k$ with
$k$ in $\{(n-3)^{(1)},\dots,1^{(1)},1^{(0)},\dots,(n-2)^{(0)}\}$.  
\end{lemma}

\begin{proof} If $I$ has such a root, then the $j$-th edge which is 
added will be an edge joining $n-1$ to such a $k$.  (Since $j\leq n-2$,
at the $j$-th step, at least one of the vertices in 
$\{(n-3)^{(1)},\dots,(n-2)^{(0)}\}$ will be available.)

On the other hand, if $\psi_{D_n}(I)$ contains such an edge with $k=k^{(0)}$,
the only possibility is that there was a root $(j,n-1)$ in $I$.  If 
$k=k^{(1)}$ then an edge from $k$ could have been added at the $k$-th
step, but this edge would not have been joining $k^{(1)}$ to $(n-1)^{(0)}$
as there would have been an available vertex with a smaller label.  
\end{proof}

We say that a vertex is the clockwise end of an edge if the vertex is
not degree zero, and the vertex to which it is attached
is closer to it counter-clockwise than clockwise.  

\begin{lemma}\label{schange} $s(I)=s(\Pan(I))$ iff $(n-1)^{(0)}$ is on the
clockwise end of an edge in $\psi_{D_n/\delta}(I)$.  
\end{lemma}

\begin{proof} 
It follows from the previous lemma that 
$s(I)=s(\Pan(I))$ iff $(n-1)^{(0)}$ is
attached to some $k$ in $\{(n-3)^{(1)},\dots,(n-2)^{(0)}\}$ in 
$\psi_{A_{2n-3}}(\widehat I)$.

Suppose  
$(n-1)^{(0)}$ is attached to some $k$ in $\{(n-3)^{(1)},\dots,(n-2)^{(0)}\}$ 
in $\psi_{A_{2n-3}}(\widehat I)$.  Observe that $(n-1)^{(0)}$ cannot be degree
zero in $\psi_{D_n/\delta}(I)$, because the edge from $(n-1)^{(0)}$ to $k$ 
is entirely within $H$.  Therefore $(n-1)^{(0)}$ is on the clockwise end
of its edge.  

Conversely, if $(n-1)^{(0)}$ is on the clockwise end of an edge in 
$\psi_{D_n/\delta}(I)$, either
it is attached to $k$ in $\{(n-3)^{(1)},\dots,(n-2)^{(0)}\}$, or else it 
is attached to $(n-1)^{(1)}$.  In fact, though, it cannot be attached 
to $(n-1)^{(1)}$ in $\psi_{D_n/\delta}$.  If it were the case that 
$(n-1)^{(0)}$ and $(n-1)^{(1)}$ were attached in  $\psi_{A_{2n-3}}(\widehat I)$, 
this edge would have been removed in $\psi_{D_n/\delta}(I)$.  
Thus $s(I)=s(\Pan(I))$. 
\end{proof}

We are now ready to prove the following result:

\begin{lemma} For $I \in \NN(D_n)$, we have that $\psi_{D_n}(\Pan(I))=
\rot(\psi_{D_n}(I))$.  
\end{lemma}

\begin{proof} By Proposition \ref{modsigmacomp}, we know that 
$\psi_{D_n/\delta}(\Pan(I))=\rot(\psi_{D_n/\delta}(I))$.  If $I$ lies in a singleton
$\delta$-orbit, this is sufficient.

Now suppose $I$ lies in a doubleton $\delta$-orbit. By Proposition 
\ref{modsigmacomp}, we know that $\rot(\psi_{D_n}(I))$ and 
$\psi_{D_n/\delta}(\Pan(I))$ differ, if at all, only in the way that the
singleton vertices are connected.  

Let $v$ be a singleton vertex in $\psi_{D_n/\delta}(I)$.  We know that
$\rot(v)$ is a singleton vertex in $\Pan(I)$.  In 
$\psi_{D_n}(I)$, suppose that $v$ is connected to $i$.  We then see
that $\rot(v)$ is connected to $i+1$, since the last two terms
in the formula cancel each other out by Lemma \ref{schange}. 
\end{proof}

\subsection{Exceptional types}
  As for noncrossing partitions in Section~\ref{sec:exceptionaltypes}, the exceptional types -- as we consider only crystallographic reflection groups, this includes for now the dihedral group $G_2$ -- were verified using a computer.

\section{Parabolic induction in the classical types}
In this section, we define the notion of parabolic induction 
for a collection of maps from $\NN(W)$ 
to $\mathcal T_W$, for $W$ a reflection group of classical type, 
and we show that
the previously defined bijections $\psi_W$ satisfy this notion of 
parabolic induction.  Further, we show that they are uniquely characterized
by this property together with their compatibility with Panyushev 
complementation and rotation.  

\subsection{Type $\An$}
First, consider the case of $W=\An$.  Pick $i$, with $1\leq i \leq n-1$.  
Removing the node $i$ from the Dynkin diagram, we obtain two Dynkin diagrams,
of types $A_{i-1}$ and $A_{n-1-i}$.  Given noncrossing handshake configurations 
$U\in \mathcal T_{A_{i-1}}$ and $V\in \mathcal T_{A_{n-1-i}}$, 
we can assemble them into a single noncrossing handshake configuration $U\ast V$ of type $A_{n-1}$, 
by adding $i$ to the labels of the vertices of 
$V$.  (In order for this
to work if $i=1$ or $i=n-1$, we define the unique 
noncrossing handshake configuration associated to
type $A_0$ to consist of two vertices, numbered $1^{(0)}$ and $1^{(1)}$, 
connected by an edge.)

Suppose that $I \in \NN(A_{n-1})$ does not have $\alpha_i$ in its support.  
We can
then write $I$ as a union of $I_1$ supported over a subset of 
$\alpha_1,\dots,\alpha_{i-1}$, and $I_2$ supported over a subset of 
$\alpha_{i+1},\dots,\alpha_{n-1}$.  

We say that a collection of maps 
$F_{A_{n-1}} : \NN(A_{n-1}) \longrightarrow \mathcal T_{A_{n-1}}$
satisfies parabolic induction if, whenever $I \in \NN(A_{n-1})$ satisfies that
the simple root $\alpha_i$ is not in the support of $I$, then 
$$F_{A_{n-1}}(I)=F_{A_{i-1}}(I_{1}) \ast F_{A_{n-1-i}}(I_{2})$$

\begin{proposition} The maps $\psi_{A_{n-1}}$ satisfy parabolic induction.
\end{proposition}

\begin{proof} This is an immediate corollary of 
Proposition~\ref{prop:supporthandshake}. \end{proof}

\subsection{Type $C_n$}
Similarly, if we remove a simple root $\alpha_i$ from a $C_n$ Dynkin diagram, 
we obtain a diagram of type $A_{i-1}$ and one of type $C_{n-i}$.  
For convenience, we use $C_1$ as a pseudonym for $A_1$ here.  In particular,
the noncrossing handshake configurations of type $C_1$ are just the 
noncrossing handshake configurations of type $A_1$.  By convention, the empty diagram is the 
unique noncrossing handshake configuration of type $C_0$.

Given
a noncrossing handshake configuration of type $U\in \mathcal T_{A_{i-1}}$ and 
$V\in\mathcal T_{C_{n-i}}$, define $U\ast V$ to consist of: \begin{itemize}
\item $U$,
\item $V$ with its labels increased by $i$,
\item $U$ with each label $j$ replaced by $2n+1-j$, and superscripts 
$(0)$ and $(1)$ interchanged.
\end{itemize}

Again, if $I \in \NN(C_n)$ and $\alpha_i$ is not in the support of $I$, we
can divide $I$ into antichains $I_1$ and $I_2$.
A collection of maps $F_{W}:\NN(W)\rightarrow \mathcal T_W$ for
$W$ of type $A$ or $C$ is said to satisfy parabolic induction if the
collection $F_{A_n}$ satisfies type $A$ parabolic induction and for 
$I \in \NN(C_n)$, whenever
$\alpha_i$ is not in the support of $I$, we have
$$F_{C_n}(I)= F_{A_{i-1}}(I_1)\ast F_{C_{n-i}}(I_2).$$

We have the following corollary of the previous proposition:
\begin{corollary}
The maps $\psi_{A_{n}}, \psi_{C_n}$ satisfy parabolic induction.
\end{corollary}

\subsection{Type $D_n$} If we remove a simple root $\alpha_i$ 
from a Dynkin diagram
of type $D_n$, for $i\ne n-1, n$ (the two antennae), then we obtain
a Dynkin diagram of type $A_{i-1}$ and a Dynkin diagram of type 
$D_{n-i}$.  Given two noncrossing handshake configurations $U\in \mathcal T_{A_{i-1}}$ and
$V\in \mathcal T_{D_{n-i}}$, we write $U\ast V$ for the diagram consisting of:
\begin{itemize}
\item The diagram $U$,
\item The diagram $V$ with its labels increased by $i$ (including the
central ones, where the increase is taken modulo 4),
\item The diagram $U$ with label $j$ replaced by $2n-1-j$, and the 
superscripts $(0)$ and $(1)$ interchanged.  
\end{itemize}
(We let $D_2$ refer to the reducible root system consisting of
two orthogonal simple roots and their negatives, and let $D_3=A_3$.  We interpret ``noncrossing handshake configuration of type $D_n$'' for $n=2,3$, using the type $D$ definition of
noncrossing handshake configuration.)

If we remove a simple root $\alpha_i$ from a Dynkin diagram of type $D_n$,
where $i=n-1$ or $n$, then we obtain a Dynkin diagram of type $A_{n-1}$.
We will define a pair of maps 
$\Ind_i:\mathcal T_{A_{n-1}}\rightarrow \mathcal T_{D_n}$, as follows.  
 
$\Ind_n(U)$ is defined to consist of the type $A$ diagram, with vertices 
$n^{(0)}$ and $n^{(1)}$ moved to the center and renamed $n$ and $n+1$, together
with the 180 degree rotation of this diagram.  
This is a type $D_n$ noncrossing handshake configuration by Lemma \ref{restatement}.  

$\Ind_{n-1}(U)$ is obtained by adding 2 to each of the labels of the central
vertices of $\Ind_n(U)$.

Again, if $1\leq i \leq n-2$, and $I \in \NN(D_n)$ does not have $\alpha_i$ 
in its support, we can define $I_1\in \NN(A_{i-1})$ and $I_2\in\NN(D_{n-i})$.
If $i=n-1,n$, and $I$ does not have $\alpha_i$ in its support, we can
simply view $I$ as an antichain in $\NN(A_{n-1})$.
A collection of maps $F_W : \NN(W) \longrightarrow \T_W$ for $W=A_n,D_n$ is said to satisfy 
parabolic induction if the collection 
$F_{A_n}$ satisfies type $A$ parabolic induction, and:
\begin{itemize}
\item[(i)] for $1\leq i \leq n-2$, if $I \in \NN(D_n)$ does not have 
$\alpha_i$ in its support, then 
$$F_{D_n}(I)=F_{A_{i-1}}(I_1)\ast 
F_{D_{n-i}}(I_2),$$ and
\item[(ii)] for $i=n-1,n$, if $I \in \NN(D_n)$ does not have $\alpha_i$ in its 
support, then 
$$F_{D_n}(I)=\Ind_i(F_{A_{n-1}}(I)).$$
\end{itemize}

\begin{proposition}
The maps $\psi_{D_n}$, $\psi_{A_n}$ satisfy parabolic induction.  
\end{proposition}

\begin{proof}
Condition (i) follows as in the previous cases. 
For condition (ii), we divide into cases.  

{\it $I \in \NN(D_n)$ has neither
$\alpha_n$ nor $\alpha_{n-1}$ in its support.}
In this case, $\overline I$ does not intersect $R$.  The result in
this case follows as in type $C_{n-1}$.  
 
{\it $I \in \NN(D_n)$ has exactly one of $\alpha_n, \alpha_{n-1}$ in its support}.
In this case, $\overline I \cap R$ consists of either one root 
$(n-1,n)$ or two roots $(j,n)$ and $(n-1,2n-1-j)$.  It follows that
$\widehat I \cap R$ consists of either zero roots or one root.  

In the former case, in the type $A_{2n-3}$ noncrossing handshake configuration associated to
$\widehat I$, there are no edges from vertices with labels at most $n-1$
to those with labels at least $n$.  It follows that the innermost edges
from $H$ to $H^c$ are connected to $n^{(0)}$ and to $(n-1)^{(1)}$, and thus
that in the $D_n/\delta$ noncrossing handshake configuration, $(n-1)^{(1)}$ 
is a singleton vertex. The other singleton vertex with label at most 
$n-1$, call it $a$, is the one that is connected to $(n-1)^{(1)}$ in the
type $A_{2n-3}$ noncrossing handshake configuration.  
Now, suppose $I$ is supported over $\alpha_{n-1}$, so $s(I)=0$.  We deduce that
$(n-1)^{(1)}$ is attached to $n-(2n-3)+2(n-1)=n+1$.  
On the other hand, if $I$ is supported over $\alpha_n$, 
$(n-1)^{(1)}$ is attached to $(n+1)+2$.  

Now consider the calculation of $\psi_{A_{n-1}}(I)$.  Up to the $n-1$-th
step, the same thing happens.  At the $n-1$-th step, there now {\it is}
an entry in the $n-1$ column (namely, $(n-1,n)$), so we mark 
$n^{(0)}$ with label $n-1$, and thus on turn $n-1$, we 
connect $n^{(0)}$ to the nearest available entry, which must be $a$,
since it and $(n-1)^{(1)}$ are the only unmatched vertices on the 
lefthand side.  On the final step, we join $n^{(1)}$ and $(n-1)^{(1)}$.  
We see that $\psi_{D_n}(I)=\Ind_{n-s(I)}(\psi_{A_{n-1}}(I))$. 

Next, consider the case that $\widehat I \cap R$ has one root in its support,
say $(i,2n+1-i)$.  
Consider the calculation of $\psi_{A_{2n-3}}(\widehat I)$ and of 
$\psi_{A_{n-1}}(I)$ in parallel.
The same thing happens in both up to the $i$-th step.  
On the $i$-th step of the $A_{n-1}$ calculation, 
the label $i$ goes onto the node $n^{(0)}$, so we connect $n^{(0)}$
to $(n-1)^{(0)}$ at this point, while for the $A_{2n-3}$ calculation, we connect
$(2n+1-i)^{(0)}$ to $(2n-i)^{(0)}$.  From here on, the calculations
run the same way up to and through the $n-1$-th step.  In both the 
calculations, there is no entry in the $n-1$-th column, so we 
connect $(n-1)^{(1)}$ to some entry on the lefthand side.
After this step, in the calculation of $\psi_{A_{2n-3}}(\widehat I)$, 
there
are two remaining unmatched vertices whose labels are at most $n-1$.  One
of them is $(n-1)^{(0)}$, while we call the other one $a$.  
It follows that
the four vertices in $H$ which will eventually be matched to vertices
in $H^c$ are, in clockwise order, the vertex attached to $(n-1)^{(1)}$,
$a$, $(n-1)^{(0)}$, and (by symmetry) $n^{(0)}$.  The two innermost edges
are therefore the ones attached to $a$ and $(n-1)^{(0)}$.  It follows that
we will connect $a$ and $(n-1)^{(0)}$ to the internal vertices, and 
$(n-1)^{(0)}$ will be connected to $n$ if $s(I)=0$ and $n+2$ if $s(I)=1$.  

On the $n$-th
step of the $\psi_{A_{n-1}}(I)$ calculation, we connect $n^{(1)}$ to the 
only available vertex, $a$.  We therefore see that 
$\Ind_{n-s(I)}(\psi_{A_{n-1}}(I))=\psi_{D_n}(I)$, as desired.  
\end{proof}

\subsection{Uniqueness of $\psi$ in the classical types}
Finally, we show that parabolic induction determines $\psi$ uniquely
in the classical cases. In this section,
we show that:
\begin{theorem} \label{th:classicalparabolicinduction}
  The only collection of bijections $F_W:\NN(W)\rightarrow \T_{W}$, for $W$ running over all classical irreducible reflection groups, that satisfy:
  \begin{itemize}
    \item[(i)] $F_W\circ \Pan=\rot \circ F_W$, and
    \item[(ii)] classical parabolic induction, as defined previously,
  \end{itemize}
  are the maps $F_W=\psi_W$.  
\end{theorem}

\begin{proof}
We have already shown that the maps $\psi_W$ do satisfy the two 
properties mentioned in the theorem; we need only show that these two
properties are sufficient to characterize these functions uniquely.  

By property (i), it suffices to know that, for any $\Pan$ orbit in 
$\NN(W)$, there is some antichain to which some parabolic induction
applies.  Expressed in those terms, it is not obvious that this is true.  
However,
thanks to the bijections $\psi_W$, it is sufficient to show that for
any $\rot$ orbit in $\T_W$, there is a noncrossing handshake configuration which could
have arisen by parabolic induction.  This is quite clear.  Let $T$ be 
a noncrossing handshake configuration of type $W$.  Pick some
edge joining two external vertices.  After applying a suitable power of 
$\rot$ to $T$, the chosen edge connects $i^{(0)}$ to $i^{(1)}$.  
In type $A_{n-1}$, this implies that $T$ comes from a parabolic induction
$A_{i-1}\ast A_1 \ast A_{n-i}$, where at most one of these is zero.  A 
completely similar approach works in type $C$ or $D$, except in the case of
$D_2$, since in that case there is a $\rot$ orbit with no edge connecting
a pair of external vertices.  However, it is easy to check that 
both the elements of that orbit arise via $\Ind$.  This completes
the proof.
\end{proof}

\section{A uniform bijection}
In this section, we prove the Main Theorem.  We will begin
with the classical types.  
Let $W$ be a reflection group of classical type, 
and let $L,R$ be a bipartition of its simple
roots.  
For each of the three classical families, we define a certain bijection
$\phi_{(L,R)}:\T_W\rightarrow \NC(W,c_Lc_R)$, which will be 
a mild variant of $\phi_W$ as defined in Section \ref{sectiontwo}.  Then
we define $\alpha_{(L,R)}:\NN(W)\rightarrow \NC(W,c_Lc_R)$ by setting
$\alpha_{(L,R)}(I)= \phi_{(L,R)}\psi_W(I)$.  We then check that this bijection
satisfies the properties demanded by the Main Theorem.

Next, we show for any reflection group, classical or not, that
a bijection satisfying the conditions
of the Main Theorem is unique, if it exists.  This completes
the proof for the classical types.  
Our uniqueness
result also gives us an explicitly computable condition to verify 
whether or not there exists a bijection satisfying the conditions of the Main Theorem for a given $W$, assuming that the bijections
are known for all parabolic subgroups.  This condition was verified by
computer for the exceptional cases, 
thus establishing the result for all types.

\subsection{Type $A_{n-1}$}
  Let $\{s_1,\ldots,s_{n-1}\}$ with $s_i = (i,i+1)$ be the generators in type $\An$, and let $c_L c_R$ be a bipartite Coxeter element. As mentioned in Remark~\ref{re:CoxeterElements}, we can cyclically label the vertices of the noncrossing handshake configurations in $\T_n$ by the Coxeter element $c_L c_R$. If $s_1 \in L$, the cyclic labelling for $\phi_{(L,R)}$ is given by \begin{align}
    2^{(0)},2^{(1)},4^{(0)},4^{(1)},\ldots,3^{(0)},3^{(1)},1^{(0)},1^{(1)}, \label{eq:sL}
  \end{align}
	and if $s_1 \in R$, the cyclic labelling for $\phi_{(L,R)}$ is given by
	\begin{align}
    1^{(1)},3^{(0)},3^{(1)},\ldots,4^{(0)},4^{(1)},2^{(0)},2^{(1)},1^{(0)}. \label{eq:sR}
  \end{align}  
  \begin{theorem}\label{th:conjectureA}
    The bijections
    \begin{align*}
      \alpha_{\An,(L,R)}: \NN(\An) &\tilde{\longrightarrow} \NC(\An,c_L c_R), \\
      \alpha_{\An,(R,L)}: \NN(\An) &\tilde{\longrightarrow} \NC(\An,c_R c_L)
    \end{align*}
    satisfy the conditions in the Main Theorem in type $A$.
  \end{theorem}
  \begin{proof}
    We will only check the first statement; the proof of the second is identical. We must check the three properties of the Main Theorem.  The initial condition is easily verified.  The $\Pan=\K$ condition follows from the facts that $\psi_{A_{n-1}}\circ \Pan = \rot\circ\psi_{A_{n-1}}$ and $\phi_{(L,R)}\circ\rot=\K\circ \phi_{(L,R)}$.   

    As we have proved the parabolic recursion for $\psi_{\An}$ in the previous section, it is left to prove the analogous statement for $\phi_{(L,R)}$. Let $T \in \T_n$ be a noncrossing handshake configuration such that $T_1 = \{i^{(0)},i^{(1)} : 1 \leq i \leq k \}$ and $T_2 = \{i^{(0)},i^{(1)} : k < i \leq n \}$ define submatchings of $T$ with vertices being labelled as in Proposition~\ref{prop:supporthandshake}. We have to show that
    $$\phi_{(L,R)}(T) =
      \begin{cases}
         \hspace{12pt} \phi_{(L_1,R_1)}(T_1) \hspace{3pt} \phi_{(L_2,R_2)}(T_2) & \text{if } s_k \in R \\
        s_k \hspace{3pt} \phi_{(L_1,R_1)}(T_1) \hspace{3pt} \phi_{(R_2,L_2)}(T_2) & \text{if } s_k \in L, \\
      \end{cases}
    $$
    where $L_{1/2} = L \cap S_{1/2}$ and $R_{1/2} = R \cap S_{1/2}$ with $S_1 = \{s_1,\ldots,s_{k-1}\}$ and $S_2 = \{s_{k+1},\ldots,s_{n-1}\}$. This results in $4$ different cases.
    \begin{itemize}
      \item[Case 1:] $s_1 \in L, s_k \in R$. In this case, the labelling is as in \eqref{eq:sL} and $k$ is even. The statement follows as the labelling of $T_1$ is given by
      $$2^{(0)},2^{(1)},\ldots,k^{(0)},k^{(1)},(k-1)^{(0)},(k-1)^{(1)},\ldots,1^{(0)},1^{(1)},$$
      and the labelling of $T_2$ is given by the remaining labels. These are exactly the labellings obtained as well for $\phi_{(L_1,R_1)}(T_1)$ and $\phi_{(L_2,R_2)}(T_2)$.
      \item[Case 2:] $s_1 \in L, s_k \in L$. In this case, the labelling is as in \eqref{eq:sL} and $k$ is odd. The labelling of $T_1$ is now given by
      $$2^{(0)},2^{(1)},\ldots,(k+1)^{(0)},k^{(1)},\ldots,1^{(0)},1^{(1)},$$
      and the labelling of $T_2$ is given by the remaining labels. It is a straightforward check that this differs from the labelling for $\phi_{(L_1,R_1)}(T_1)$ and $\phi_{(R_2,L_2)}(T_2)$ by having the labels $(k+1)^{(0)}$ and $k^{(0)}$ interchanged. This corresponds exactly to the additional factor $s_k$.
    \end{itemize}
		The remaining two cases for $s_1 \in R$ are solved in the analogous way.
  \end{proof}

\subsection{Type $C_n$}
  As above, the bipartite Coxeter elements in type $C_n$ can be obtained from bipartite Coxeter elements in type $A_{2n-1}$, where $-i$ and $2n+1-i$ are identified. The bijection in type $C$ then follows as a simple corollary from the construction in type $A$.
  \begin{corollary}
    The bijections
    \begin{align*}
      \alpha_{C_n,(L,R)}: \NN(C_n) &\tilde{\longrightarrow} \NC(C_n,c_L c_R), \\
      \alpha_{C_n,(R,L)}: \NN(C_n) &\tilde{\longrightarrow} \NC(C_n,c_R c_L)
    \end{align*}
    satisfy the conditions in the Main Theorem in type $C$.
\end{corollary}

\subsection{Type $D$}
  Exactly the same argument as in type $A_{n-1}$ 
applies to the bipartite Coxeter elements in type $D_n$. Those are obtained from the bipartite Coxeter element in type $\An$ by adding $s_n = (n-1,-n)$ to $L$ if $n$ is even and to $R$ if $n$ is odd. E.g., in type $D_4$, we obtain the cyclic labelling on the outer circle for $c_L c_R$ given by
  $$2^{(0)}, 2^{(1)}, -3^{(0)}, -3^{(1)}, -1^{(0)}, -1^{(1)}, -2^{(0)}, -2^{(1)}, 3^{(0)}, 3^{(1)}, 1^{(0)}, 1^{(1)},$$
  and the inner circle labelling by $4^{(0)},4^{(1)},-4^{(0)},-4^{(1)}$. The labellings for $c_R c_L$ are again given by reflecting the labels at the diagonal through $1^{(1)}$.
  \begin{corollary}
    The bijections
    \begin{align*}
      \alpha_{D_n,(L,R)}: \NN(D_n) &\tilde{\longrightarrow} \NC(D_n,c_L c_R), \\
      \alpha_{D_n,(R,L)}: \NN(D_n) &\tilde{\longrightarrow} \NC(D_n,c_R c_L)
    \end{align*}
    satisfy the conditions in the Main Theorem in type $D$.
  \end{corollary}
  \begin{proof}
    The proof follows the same lines as the proof in type $A$, with the additional check for the cases in which $s_{n-1}$ or $s_n$ are not contained in the support of an antichain $I \in \NN(D_n)$. Using Theorem~\ref{th:classicalparabolicinduction} in type $D_n$, this check is straightforward.
  \end{proof}

\subsection{Uniqueness}
We now establish uniqueness of the bijections satisfying the conditions of
the Main Theorem.
For $I$ which has less than full support, $\alpha_{(L,R)}(I)$ is determined by
parabolic induction. By $\Pan=\K$, there is likewise no choice for
$\alpha_{(L,R)}(J)$ for any $J$ in the $\Pan$-orbit of $I$.  We saw in the
classical types,  
in the proof of Theorem~\ref{th:classicalparabolicinduction}, that  
every $\Pan$-orbit in $\NN(W)$ contains an antichain which does not have full support. 
This fact can also easily be checked (by computer) for the exceptional types.  
Therefore, there is at most one $\alpha_{(L,R)}$ satisfying the 
conditions of the Main Theorem.

\subsection{Exceptional types}
The argument above for uniqueness, actually proves more: it essentially
gives a candidate bijection.  
Suppose that bijections as in the Main Theorem
have already been defined for all proper parabolic subgroups of $W$.  
For each $\Pan$-orbit $\mathcal O$ in $\NN(W)$, pick an antichain $I_{\mathcal O}\in \mathcal O$
which does not have full support, and define 
$\alpha_{(L,R)}(I_{\mathcal O})$ by parabolic induction.  Now extend 
the definition of $\alpha_{(L,R)}$ to all of $\mathcal O$ by $\Pan=\K$.  
We now have a candidate for a map satisfying the Main Theorem's condition and,
as in the uniqueness argument above, if there is any map satisfying
the conditions of the Main Theorem for $W$, it must be this one.  
The fact that this map really is a bijection satisfying all three
of the properties of the Main Theorem can now be verified by 
computer (and has been verified) in the exceptional types.  This completes the proof of 
the Main Theorem.

\section{A proof of the Panyushev conjectures}

In this final section of the paper, we will use combinatorial results described in the previous sections to prove the Panyushev conjectures. The first proposition follows directly from the uniform description of the bijection.
\begin{proposition}
	Part (i) of the Panyushev conjectures holds: $\P^{2h}$ is the identity map on $\NN(W)$.
\end{proposition}
\begin{proof}
	This follows from the connection to the Kreweras complementation and the fact that $\K^{2h}$ is the identity map on $\NC(W)$.
\end{proof}

For all remaining proofs, we use the combinatorics obtained for the classical types, and computer checks for the exceptionals. To prove (ii) of the Panyushev conjectures, it remains to show that $\K^h$ acts on $\NN(W)$ by the involution induced by $-\omega_0$. Thus, we have two cases, depending on how 
$-\omega_0$ acts on Dynkin diagrams:
\begin{itemize}
	\item[(iia)] $\K^h$ acts trivially on $\Phi$ in type $C_n, D_{2n}, F_4, E_7,$ and $E_8$.
	\item[(iib)] In the remaining types $\An, D_{2n+1},$ and $E_6$, the action of $\K^h$ is induced by the involution on the Dynkin diagram (called $\delta$ in types $A$ and $D$).
\end{itemize}
\begin{proof}[Proof of part (ii) of the Panyushev conjectures]
	In types $A$ and $C$, (iia) and (iib) follow from the symmetry property of noncrossing handshake configurations (see Lemma~\ref{lem:Binvolution}). In type $D$, (iia) and (iib) follow from the facts that rotating a type $D_n/\delta$ noncrossing handshake configuration by $2(n-1)$ steps yields the same configuration, but to obtain the same $D_n$ noncrossing handshake configuration, it is also necessary to ensure that the number of rotations applied yields a half-turn of the $4$ inner vertices. 
Type $E_6$ was checked with a computer.  The statements for the remaining exceptional types can be verified using the orbit lengths found in Section~\ref{sec:exceptionaltypes}.
\end{proof}
\begin{proof}[Proof of part (iii) of the Panyushev conjectures]
First we consider type $\An$.  Pick 
a noncrossing handshake configuration $X$, and consider $X, \K(X),\dots,\K^{2n-1}(X)$.  Each edge $e$ in
$X$ appears (rotated) in each of these noncrossing handshake configurations, and we see that some
endpoint of $e$ is labelled with $(0)$ and marked in $n-1$ of these noncrossing handshake configurations.  
In a given noncrossing handshake configuration, the number of  vertices labelled with $(0)$ and marked
is exactly the number of positive roots in the corresponding antichain, so we see that
the total number of positive roots in the antichains corresponding to these $2n$ noncrossing handshake 
configurations is $n-1$ times the number of edges, which is $n$.  It follows that the
average number of positive roots in the corresponding $\Pan$ orbit is $(n-1)/2$.  

  The easiest way to prove the result for type $C_n$ is the following: 
it is straightforward to check that every second antichain in a Panyushev orbit contains a positive root of the form $(i,\overline{i})$. As type $A_{2n-1}$ folds to the type $C_n$, the total number of antichains in an orbit in type $C_n$ is given by
  $$\frac{\frac{4n}{2}\frac{2n-1}{2} + 2n}{4n} = \frac{n}{2}.$$
  Here, the nominator contains $4n\frac{2n-1}{2}$ which is the orbit size (without symmetry) times the average number of elements in the orbit in type $A_{2n-1}$, the division by $2$ comes from the folding, and the correction term $2n$ comes from the centered element in every other orbit which is not folded. The $4n$ in the denominator is again the size of the orbit. (If we have a $k$-fold symmetry, all three pieces obtain a factor of $1/k$.) This completes the proof in type $C$.

  In type $D$, the situation is again a little more involved. We will work in terms of $D_n/\delta$ configurations.  
There are two different cases, based on whether or not there are four 
isolated vertices on the outside.  

Suppose first that there are not.  Each such $D_n$ antichain 
corresponds to a $C_{n-1}$ antichain, and the Panyushev map respects this
folding action.  Thus, a Panyushev orbit of such $D_n$ antichains 
corresponds to a Panyushev orbit of $C_{n-1}$ antichains; the
average number of roots present in these $C_{n-1}$ antichains is 
$(n-1)/2$.  The $D_n$ antichain $I$ corresponding to a $C_{n-1}$ antichain 
$I'$ is just the inverse image of $I'$ under the folding
map from $\Phi_{D_n}$ to $\Phi_{C_{n-1}}$.  The number of elements in 
$I$ equals the number of elements in $I'$, plus the number of elements
in $I'$ whose inverse image consists of two roots; there will be either
one or zero such roots in $I'$.  We observe that there is such a root 
in $I'$ iff $n^{(0)}$ is marked.  As we rotate 
$\psi_{C_{n-1}}(I')$ through a full rotation,
each edge of the configuration is connected to vertex $n^{(0)}$ twice, once
at each of its endpoints, and it is easy to see that once we will
have $n^{(0)}$ marked, while once it will be unmarked.  Thus, the
average effect of passing from $I'$ to $I$ is to add $\frac12$ to the
size of the antichains, resulting in an average size of $n/2$ as desired.

Now suppose that there are four isolated vertices in 
$\psi_{D_n/\delta}(I)$.  We consider first the average size of 
$\hat I$ (which, we recall, is an antichain of type $A_{2n-3}$).
  Recall that, as we consider $\psi_{A_{2n-3}}(\hat I),
\psi_{A_{2n-3}}(\widehat {\Pan(I)}),\dots$, the effect is to rotate the
noncrossing handshake configuration except that there is one pair of edges which,
at a certain point, gets switched, and then eventually switches back; in a full rotation
($4n-4$ steps) this happens twice.  

Consider first an edge which is not involved in the switching.  It
contributes a marked vertex $2n-3$ times (out of the $4n-4$ rotations).  
Now consider the pair of edges that are involved in the switching.  
One verifies directly that they contribute, together, $4n-8$ marked vertices.  
The average size of the antichains $\widehat I, \widehat{\Pan(I)},$ etc.,
is $[(2n-4)(2n-3)+(4n-8)]/(4n-4)= (4n^2-10n+4)/(4n-4)$.  

We next consider the average size of the sets $\overline I$, 
$\overline{\Pan(I)}$, etc.  Each of these contains one more root than the
corresponding antichain $\widehat I, \widehat{\Pan(I)}$, etc., so the
average size of these sets is $(4n^2 - 6n)/(4n-4)$.

Next we consider the relationship between the size of $\overline I$ and the
size of $I$.  The size of $I$ is $|\overline I|/2$, plus a correction of 
$\frac12$ if $\overline I$ has an element on the central diagonal.  
Over $4n-4$ rotations, the correction will appear $2n$ times (i.e. two more
than half the time).  
The reason for this is that, if $I$ is such that $\widehat I$ and
$\widehat {\Pan(I)}$ differ by a switch of the edges, then neither of 
them will have an element on the central diagonal.  
We see this because 
of the fact that the switching edges are the most internal among those
connecting $H$ to $H^c$ in 
$\psi_{A_{2n-3}}(\widehat I)$.  
Now $\widehat I$ has
no element on the central diagonal iff $\overline I$ does have an element
on the central diagonal.  

It follows that the number of elements in an antichain, averaged over a 
$\Pan$-orbit,  
is $(4n^2-4n)/(8n-8)=n/2$.
\end{proof}

\section*{Acknowledgements}
The authors would like to thank Christian Krattenthaler and Vic Reiner for helpful discussions.

The Main Theorem was verified in the exceptional types with the help of J.~Stembridge's {\tt coxeter} and {\tt posets} packages for {\tt Maple}.

The orbit sizes for the Kreweras and the Panyushev complements were calculated using {\tt SAGE} \cite{sage}.

During the time that he worked on this paper, D.A. was supported by NSF Postdoctoral Fellowship DMS-0603567 and NSF grant DMS-1001825.

C.S. was supported by a CRM-ISM postdoctoral fellowship. He would like to thank the Fields Institute for its hospitality during the time he was working on this paper.

H.T. was supported by an NSERC Discovery Grant.  He would like to thank the Norges teknisk-naturvitenskapelige universitet and the Fields Institute for their hospitality during the time he was working on this paper.

  \bibliographystyle{amsplain}
  \bibliography{../../bibliography}
  
\end{document}